\documentclass[final,3p]{elsarticle}
\newcommand{\fref}[1]{(\ref{eq:#1})}
\usepackage[utf8]{inputenc}
\usepackage{bbm}
\usepackage{bbold}
\usepackage{amsmath}
\usepackage{indentfirst}
\usepackage{graphicx}
\usepackage{subfigure}
\usepackage{xcolor}
\usepackage{fancyhdr}
\usepackage{lastpage}
\usepackage{calc}
\usepackage{times}
\usepackage{lipsum}
\usepackage{wrapfig}
\usepackage{setspace}
\usepackage{sectsty}
\usepackage{hyperref}
\usepackage{longtable}  
\usepackage{epstopdf} 

\usepackage{adjustbox}
\usepackage{array}
\usepackage{booktabs}
\usepackage{multirow}

\newcolumntype{R}[2]{%
    >{\adjustbox{angle=#1,lap=\width-(#2)}\bgroup}%
    l%
    <{\egroup}%
}

\usepackage{amssymb}

\journal{Solar Energy}

\begin{document}

\begin{frontmatter}




\title{Optimal Switchable Load Sizing and Scheduling for Standalone Renewable Energy Systems }


 \author[rvt]{Abdulelah H. Habib\corref{cor1}}
 \ead{ahhabib@ucsd.edu}
 
 \author[rvt]{Vahid R. Disfani}
 \ead{disfani@ucsd.edu}
 
 \author[rvt]{Jan Kleissl }
 \ead{jkleissl@ucsd.edu}
 
  \author[rvt]{ Raymond A. de Callafon }
 \ead{callafon@ucsd.edu}
 
 \cortext[cor1]{Corresponding author: Abdulelah Habib, Tel.: +1 (858) 353 4278, Fax: +1 (858) 534-7078.}

  \address[rvt]{Department of Mechanical and Aerospace Engineering (MAE), University of California, San Diego, 9500 Gilman Drive MC 0411,
La Jolla, CA 92093-0411}

\begin{abstract}
The variability of solar energy in off-grid systems dictates the sizing of energy storage systems along with the sizing and scheduling of loads present in the off-grid system. Unfortunately, energy storage may be costly, while frequent switching of loads in the absence of an energy storage system causes wear and tear and should be avoided. Yet, the amount of solar energy utilized should be maximized and the problem of finding the optimal static load size of a finite number of discrete electric loads on the basis of a load response optimization is considered in this paper. The objective of the optimization is to maximize solar energy utilization without the need for costly energy storage systems in an off-grid system. Conceptual and real data for solar photovoltaic power production is provided the input to the off-grid system. Given the number of units, the following analytical solutions and computational algorithms are proposed to compute the optimal load size of each unit: mixed-integer linear programming and constrained least squares. Based on the available solar power profile, the algorithms select the optimal on/off switch times and maximize solar energy utilization by computing the optimal static load sizes. The effectiveness of the algorithms is compared using one year of solar power data from San Diego, California and Thuwal, Saudi Arabia. It is shown that  the annual system solar energy utilization is optimized to 73\% when using two loads and can be boosted up to 98\% using a six load configuration.

\end{abstract}


\end{frontmatter}
\section*{Nomenclature}

\begin{tabular}{ll}
$S(t_k)$ & Solar power data sampled at $t_k$ and normalized by the maximum solar generation $S_{max}$\\
$S_{max}$ & Maximum solar generation\\
$\Delta _T$ &is the sampling rate \\
T & is the length of the solar power data\\
k& is an array of size $T$ \\
$x^i$ &static load size unit\\
$n$ &is the number of loads\\
$E(t_k)$ &power mismatch or difference of solar power $S(t_k)$ and power used by the loads $u(t_k)  x$\\
$u^i_k$&binary number represent the switching of the loads\\
$N$&Number of possible combination of units \\
$U$ &is the permutation matrix of size $[n,N]$\\ 
\end{tabular}
\section{Introduction}
\label{sec:intro}

Increasing global energy demand and human population growth have triggered a need for standalone renewable applications. Recent estimates show that 1.4 billion people do not have access to energy services and one billion are suffering from unreliable electricity services \cite{IEA}. Standalone application of clean energy, (E.g., fresh water pumping), has become more critical for humanity \cite{IEA,UN}. Often such systems are powered by solar photovoltaic (PV) due to ubiquitous high solar resource availability and scalability. However, solar production exhibits high variability over a broad range of time scales \cite{nrelvar}. Power variability is the main obstacle facing solar energy in standalone or {islanded} mode applications. {High penetrations} of solar power sources create {large power swings which influence} electric power quality \cite{pvvoltage}, and can cause loss of load or generation curtailment \cite{windsolarchalng}.   Variability of solar PV generation is a result of seasonal and diurnal changes in the sunpath as well as short-lived cloud cover. Solar variability limits the operation of off-grid loads at maximum capacity \cite{griduncer,reviewstandalonePV,sreeraj2010design}. 

Optimal load switching can be applied to microgrids with any hybrid forms of renewable energy resources such as solar and wind \cite{microgirdoptim,indiasizing,micoopti1,micoopti2} to capture as much renewable energy as possible. Although partial or modulated load operation is conducive to the problem, there are numerous types of load units which can only be switched on or off, such as non-dimmable lighting, standard electric motors, and Magnetic Resonance Imaging (MRI) machines at hospitals and load aggregation such as demand side management \cite{loadagg,hotwater}. Dispatching such binary load units, which are referred as switchable loads hereafter, to follow available renewable energy resources have been discussed in the literature for different microgrid applications such as water desalination \cite{desal}, pumping systems \cite{PVpump}, irrigation systems \cite{pvwater}, and cooking appliances\cite{standalone1,standalone2,PVstandalone}. 

Different optimization techniques are used for planning and design of such systems. For instance, mixed-integer linear programming (MILP) has been used in many fields, such as unit commitment of power production \cite{MILPpower} and power transmission network expansion \cite{MIPtran,MILPpowertrans}, as well as scheduling problem of the generation units in off-grid in order to maximize supply performance of the system \cite{MNLPmicrogrid}. Nonlinear approaches have also been applied to load scheduling \cite{Nonlinearsch}. For example, neural networks and genetic algorithms have been applied to size stand-alone PV \cite{PVstandalone2,artifPV,gaPV}. The on/off control optimization problem is similar to the unit commitment problem in power systems and bio-fuel\cite{onoff,onoff2}.
However, limited studies have been conducted on the optimal load sizing in a standalone (islanded) grid application with switchable loads. 
Most of other research has been in the demand/supply side while very few looked into the unit/load sizing for many reasons, such as, the load is assumed to be fixed and has to meet by any supply way \cite{Ashok20071155} or the accessibility of designing load for certain application is harder and not easy process.
This work focuses on optimal load sizing for standalone applications in rural areas or off-grid sites. While the present paper assumes an off-grid system, similar challenges exist for a power system with a weak grid connection, i.e. with a line carrying capacity that could only balance variability that is a small fraction of local solar generation or load capacity.

{Energy storage systems (ESS) have been applied to solve the variability  challenges \cite{pickard2012addressing,storagereview}.An alternative or complementary approach is optimal sizing and scheduling of load units which follow power generation variability to maximize solar energy utilization and load uptime. Clearly, the solar energy utilization could be improved with an ESS, but an ESS that eliminates solar variability would need to be large enough to store several days' worth of solar power which is uneconomical at present. Smaller ESS would experience significant cycling and deep discharge events if not properly maintained, increasing maintenance costs and requiring replacement much before the end-of-life of a PV system}. Our objective is to improve solar utilization without an ESS and use load demand response only and show that high efficiencies can still be obtained. In practice, a combination of a small ESS with high cycle life such as an ultracapacitor ESS and the proposed load sizing and scheduling system would probably be the best solution. The ESS would absorb solar variability at time scales of seconds to minutes while the loads would balance variability at longer time scales. This approach would allow limiting ESS energy capacity making it more economical. While practical challenges of implementing such a system are significant, e.g. in maintaining system stability during switching, this paper focuses on the critical algorithmic work that permits such a system to operate efficiently and economically.

In a properly planned system the solar system would be optimally sized to power the load required for the intended application. This paper does not consider this scenario. Often in practice the conditions are not as plannable. Load growth will occur and a solar power system may be initially oversized to accommodate such growth. Sizing the solar system may also be limited by land ownership and topographic constraints. The solutions proposed in this paper apply in such a context where solar capacity is fixed and loads are sized to optimize solar energy utilization.

This paper proposes an optimization model to capture the maximum amount of variable solar generation, which sizes and schedules a finite number of loads to track available solar PV power. The objective is to maximize solar utilization, given the projected power generation of the renewable energy resources. Here, solar utilization is defined as percentage of energy captured by the units over total solar energy produced. This is akin to terms such as solar utilization factor \cite{solarutlizaton} and loss of power supply (LPS) \cite{LPSP} which are commonly used in the literature. The loads are assumed to switch between a binary "on" or "off" statuses, where both the switching times and the size of the static power demand (static load size) determines the ability to track available solar power. 

The main contribution of this paper is to develop both analytical solutions and computational approaches based on Equality Constrained Least Squares (ECLS), Inequality Constrained Least Squares (ICLS), and Mixed-Integer Linear Programming (MILP) in order to solve the optimization problem. The rest of the paper is organized as follows. The mathematical formulation is given in Section \ref{sec:problem} along with an analytical example and the motivation for a computational procedure for optimal load size selection. One year of solar resource data for San Diego is analyzed and discussed in Section\ref{sec:solardata}. 
Section \ref{sec:optapp} presents different computational procedures for optimal load size selection based on a bi-linear optimization problem involving a mix of binary and real numbers. The simulation results are presented and discussed in Section ~\ref{sec: Results}. Finally, Section \ref{conc} concludes the paper.

\section{Problem Formulation}
\label{sec:problem}
The sizing and scheduling problem computes the distribution of the optimal load size of a finite number load units, given an available (solar) power profile. For the optimal load size selection, the loads are assumed to operate in a binary manner, off or on, and therefore only the static load size is optimized. 

\subsection{Static Load Response Optimization Problem}
To formalize the notation for the optimization approach presented in this paper, we assume knowledge of the (solar) power delivery $S(t)$ sampled at regular time intervals $t=t_k = k \Delta_T$, where $1/\Delta_T$ is a fixed sampling frequency and $k$ is the sample index. In this way, we have a data set of $T$ points on the solar power production $S(t_k),~ k=1,...,T$. Typically, $S(t_k)$ is close to a daily periodic function and $S(t_k)\geq0$ over a daily time interval $t_k \in [t_b,t_e]$, where $b$ is the beginning and $e$ is the ending of the day, with a maximum value
\[
S_{max} = \max_k S(t_k),
\]
that is typically equal to the AC rating of the solar power system. {The data} $S(t_k)$ and $S_{max}$ may be available from historical solar power measurements or from solar radiation measurements in conjunction with a model of the solar power conversion efficiency.  

In the static load response optimization we consider $n$ loads, where each load $i=1,...,n$ is characterized only by a {static} power value $x^i$ that can be either turned on or off. Given the number $n$ of loads, the objective of the static load response optimization is to find the optimal distribution of static load values $x^i,~ i=1,2,\ldots,n$ so that the time sampled solar power delivery $S(t_k)$ can be approximated as closely as possible to maximize the energy captured.

\begin {figure}[ht]
\graphicspath{ Plots/ }
\includegraphics[width=.4\columnwidth]{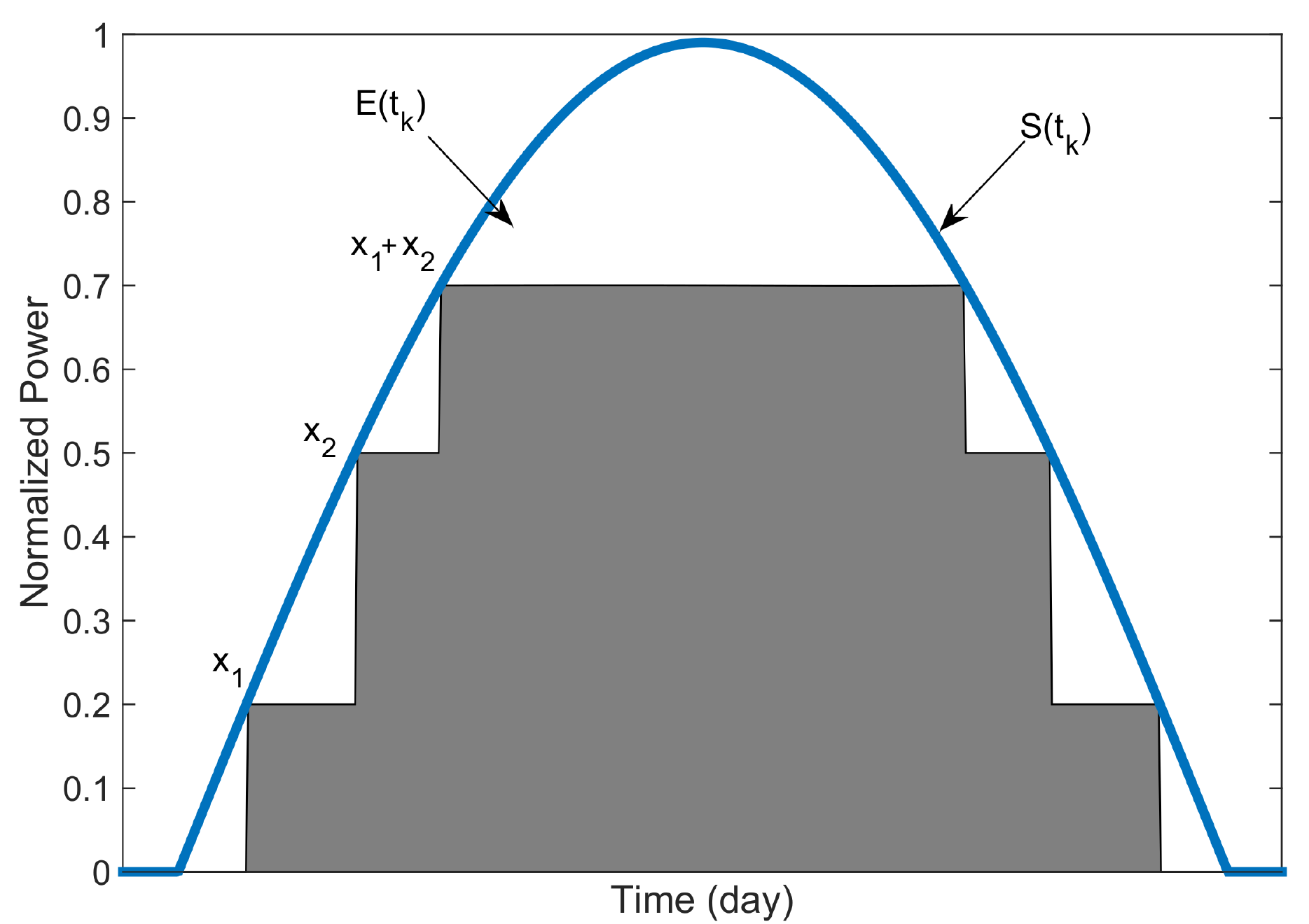}
\centering
\caption{Problem illustration for one clear and symmetric day. Loads $x^i$ are scheduled hourly to follow the increase or decrease in solar power generation $S(t_k)$.}
\label{figure:min}
\end {figure}

As indicated in Figure~\ref{figure:min}, the power mismatch $E(t_k)$ at any time 
$t_k$  can be characterized by
\[
E(t_k)= S(t_k)-\sum_{i=1}^{n} u^i_kx^i
\]
where $u_k^i \in [0,1]$ are $n$ binary numbers reflecting the on/off {switch state} of the individual loads $i=1,\ldots,n$ with their {(to be determined) static load size $x^i$}. Defining the vectors
\begin{equation}
\begin{array}{rcll}
u_k &=& \left [ \begin{array}{cccc} u_k^1 & u_k^2 & ... & u_k^n \end{array} \right ],~ & {u_k^i} \in [0,1]\\ x &=& \left [ \begin{array}{cccc} x^1 & x^2 & ... & x^n \end{array} \right ]^T,~ & x^i > 0
\end{array}
\label{eq:vars}
\end{equation}
the static power mismatch $E(t_k)$ at a particular time $t_k$ can be written with an inner product
\[
E(t_k) = S(t_k) - u_k x,
\]
of the {time dependent binary switch state vector $u_k$ and the static load size distribution vector $x$}. Static load response optimization can now be written as
\begin{equation}
\mbox{arg} \min_{u_k,x} \sum_{k=1}^N E(t_k)^2,~~ E(t_k) = S(t_k) - u_k x 
\label{eq:opt}
\end{equation}
where the variables $u_k$ and $x$ are given in \fref{vars}.
The optimization in \fref{opt} is a least squares optimization in which both the {time dependent binary switch state vector $u_k$ and the static load size distribution vector $x$} must be determined on the basis of the $T$ data points on the solar production $S(t_k),~ k=1,\ldots,T$.

Clearly, the least squares optimization in \fref{opt} is non-standard for several reasons. First of all, the error $E(t_k)$ is bi-linear due to the product of the optimization variable {$u_k$} and $x$. Furthermore, the optimization variable $u_k$ is a binary vector, whereas {the elements $x^i$ of the static load size distribution vector $x$ must likely satisfy (linear) constraints
\begin{equation}
x^i\geq x^{i-1} \geq 0
\label{eq:ordering}
\end{equation}
to ensure a unique load distribution solution with real valued positive loads}. An additional linear constraint
\begin{equation}
\sum_{i=1}^n x^i = |x|_1 = {C} S_{max}
\label{eq:cons}
\end{equation}
where $C$ can be chosen in the range of ${0.5<C<1}$ ensures that the sum of the load distribution is bounded to avoid oversizing of the loads in trying to match the anticipated maximum power production $S_{max}$. Finally, the number of loads ($n$) also needs to be determined. It is clear that a larger value {$n$} of will enable smaller power mismatch errors $E(t_k)$ but would likely increase the investment cost as the per kWh cost of a load unit typically decreases with the size of the unit.

{The process can be described in a flow chart as Figure \ref{figure:flowchart} illustrates, starting with annual solar irradiance data for at least one year to capture the seasonal changes. Then followed by a power model to compute solar power and then sort solar power annual data followed by inputting number of desired units then the optimal unit sizing will be computed then followed by a daily unit scheduling.
\begin{figure}[ht]
\graphicspath{ {Plots/} }
\includegraphics[width=.3\columnwidth]{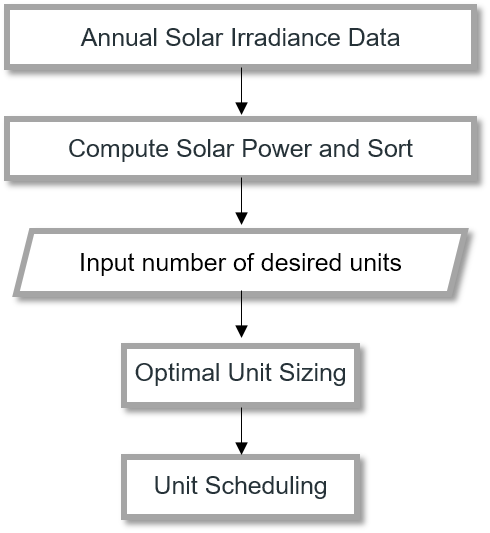}
\centering
\caption{{Flow chart of the proposed optimization process in this paper.}}
\label{figure:flowchart}
\end{figure}
}

\section{Solar Resource Data and its Distribution}
\label{sec:solardata}
To capture the interannual variability of solar irradiance, solar resource data should be collected for several years, such as in the production of a typical meteorological year (TMY). In this paper, only one year of solar power generation is used to demonstrate the model application, but most large solar system developers rely on multidecadal modeled power production based on site adoption of long-term satellite records with short-term local measurements for their financial calculations \cite{longtermpred}. Such long-term data would be preferable in practice although interannual variability of solar energy generation is small. For example, \cite{GermanSolar} specifies that the interannual variability of GHI for 7-10 years of measurement at Potsdam, Germany and Eugene, USA is about 5\%. 

\begin {figure}[ht]
\graphicspath{ {Plots/} }
\includegraphics[width=.5\columnwidth]{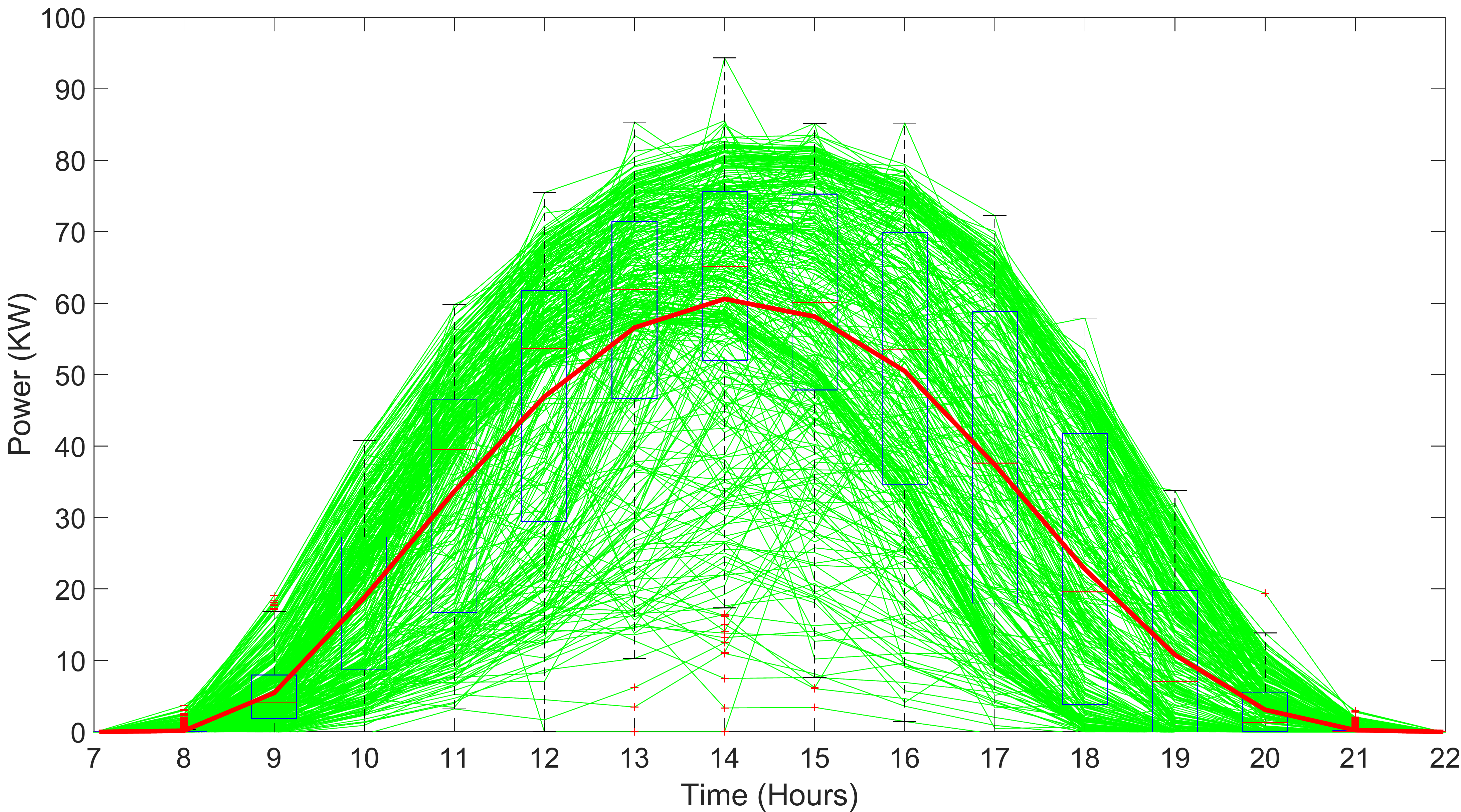}
\centering
\caption{365 days of solar PV output data from the UC San Diego campus (green) superimposed with the annual average (red). Hourly boxplots (black) show the median, $25^{th}$ and $75^{th}$ percentiles, and range.}
\label{figure:boxplot}
\end {figure}

Here, only one year of data was available which allows characterizing most of the important seasonal and diurnal variability.
{\subsection{Case Study 1: San Diego, USA}}
Multidecadal PV power projections are typically based on the measured and modeled GHI which is transposed to the direct, diffuse, and reflected radiation at the plane-of-array\cite{pvmodel1,pvmodel2,pvmodel3} and input into a PV performance model. We bypass the complexity in the PV power modeling by using a solar power generation dataset available at two sites. AC power production and GHI data were collected from a (91.6~kW$_{\rm DC}$ and 100~kW$_{\rm AC}$) fixed tilt (non-tracking) polycrystalline PV system (Figure \ref{figure:boxplot}). The system was installed at 10 $^{\circ}$ tilt and facing south, at the UC San Diego campus at $32^{\circ}53'01.4"$N $117^{\circ}14'22.6"$W. Data for one year, {\it{i.e.,}} May 2011 through April 2012, were used and averaged over 15 minutes. The raw data is available at 1~s resolution and the algorithm can be applied to data at any temporal resolution. Nevertheless, switching loads over such short timescales is generally impractical and we assume instead that a small energy storage system modulates high frequency solar variability to create a supply that is stable over 15~min intervals. The solar power data for {these specific sites are} not symmetric over a day, and overcast conditions occur more frequently in the mornings. 

\begin {figure}[ht]
\graphicspath{ {Plots/} }
\includegraphics[width=9cm]{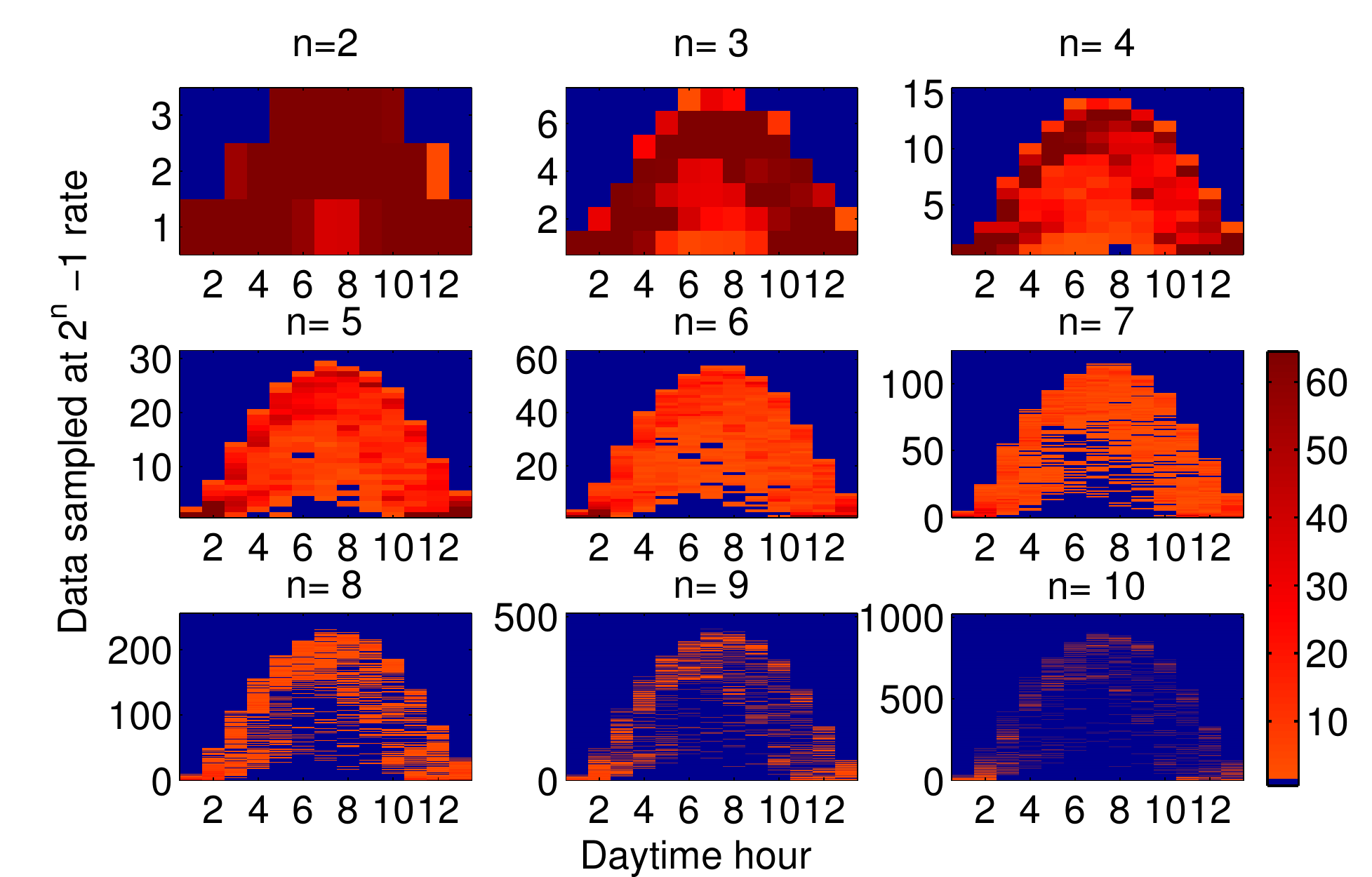}
\centering
\caption{2D histogram of number of occurrences of a combination of discrete loads for $n=2$ to 10 as a function of time of day for one year of solar power generation data shown in Figure \ref{figure:boxplot}. Only hours with non-zero power data, i.e., only daytime hours, are shown. For example, for $n=2$, either two units (at different capacity) can be turned on individually or together resulting in 3 discrete loads. For $n=10$, almost 1,000 discrete load levels exist, The colors show how often each combination of loads is run to utilize the energy from solar PV generation. }
\label{figure:Hist}
\end {figure}

The impacts of solar PV generation variability on the utilization of different combinations of discrete loads is illustrated through 2-dimensional (2D) histograms in Figure \ref{figure:Hist}. In each subplot, similar to Figure \ref{figure:boxplot}, one year of solar power data is superimposed over one day and different numbers of units $n$ are used to track the solar power generation. 
The quantization idea of the solar power data $S(t_k)$ for different numbers of  units is illustrated. $N = $ $2^n-1$ is the number of discrete combinations of units for each case described in \fref{vars}, which implies that with more units more discrete load levels can be served resulting in a better match with the solar generation. Qualitatively, the best combination of load sizes is expected to be the one that is able to best track the solar generation levels (or minimize power mismatch as in \fref{opt}) that are (i) large and (ii) occur frequently.
 
{\subsection{Case Study 2: Thuwal, Saudi Arabia}
Another case was selected to prove the robustness of the algorithms to different input data. Thuwal, north of Jeddah, is a city located in the west coast of Saudi Arabia as shown in Figure 2 of  \cite{habibasme}. Thuwal solar meteorology is predominantly clear and at a lower latitude and is therefore quite different from San Diego with days such as clear day (D1) in Figure \ref{figure:results} being more common. Data from a monocrystalline Silicon solar PV power plant at ($22^{\circ}18'28.5"$N $39^{\circ}06'17.1"$E) with tilt 20$^{\circ}$ and azimuth of 133$^{\circ}$ and 145$^{\circ}$ (split in two different arrays) was collected by King Abdullah University of Science and Technology (KAUST). }

\section{Optimization Techniques}
\label{sec:optapp}
{This section is initiated by an analytical solution for a small number $n\leq 2$ of loads is presented when the available power follows a symmetric daily profile. The analytical approach is followed by the motivation for a computational procedure to compute optimal load size distribution for a larger number $n \geq 2$ of loads when the available daily power profile is non-symmetric. After that three different optimization techniques are presented, discussed and compared. }
\subsection{{Analytical and Motivating Example}}
\label{sec:motivation}

\begin{figure}[ht]
\graphicspath{ {Plots/} }
\includegraphics[width=.5\columnwidth]{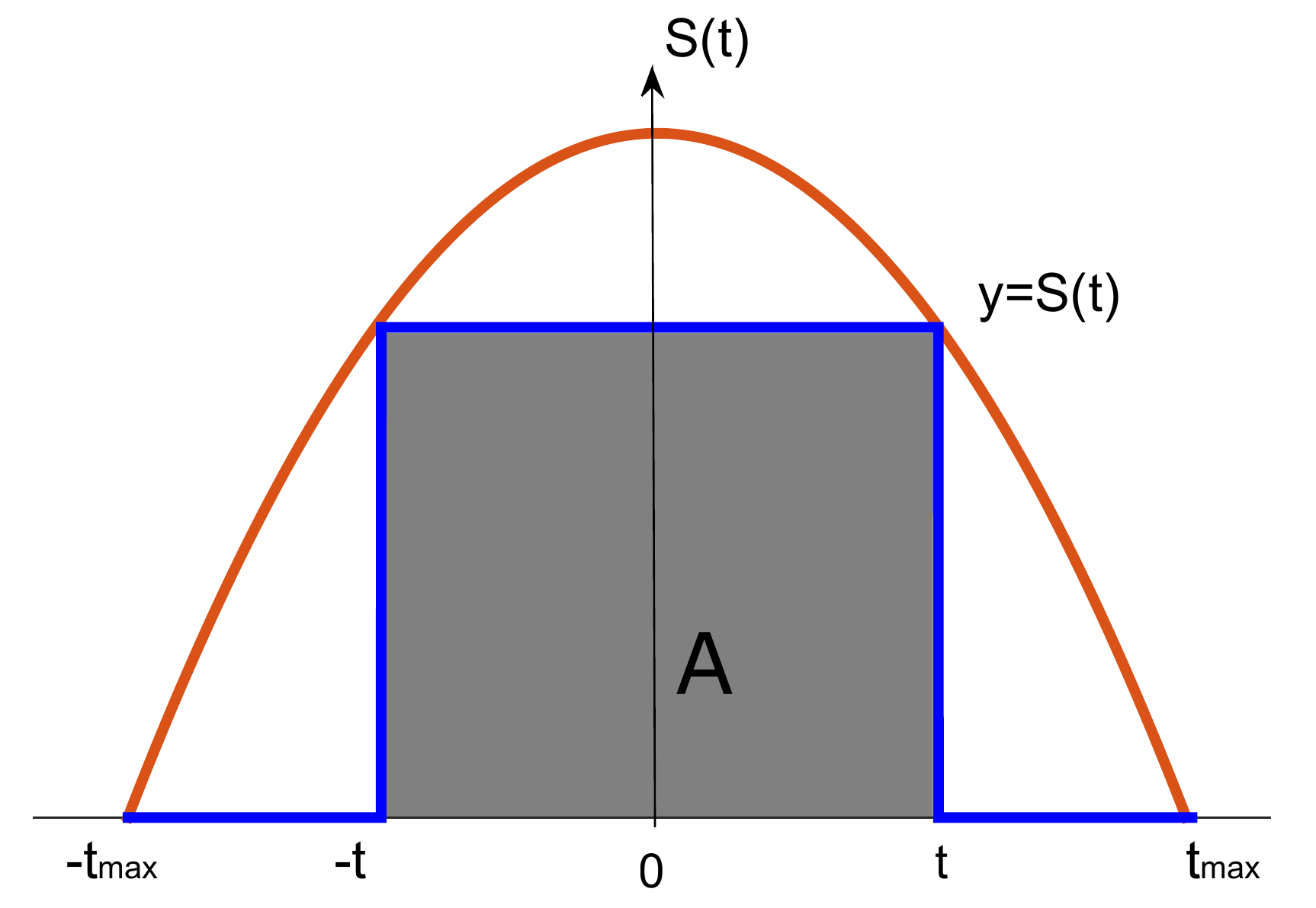}
\centering
\caption{Single static unit {sizing} and scheduling optimization over a symmetric power function $S(t)$, resembling solar power output on a clear day.}
\label{figure:con1}
\end {figure}

Consider a (symmetric) time dependent power function $y=S(t)$ which must be followed by the rectangular power demands created by a simple on/off switching of a static load over a specified time period. System efficiency is optimized by finding the largest rectangular window (representing energy demand by the switchable load units) to be drawn under the power function $S(t)$. For a single load subjected to a symmetric $S(t)$, this problem reduces to selecting an optimal on/off switch time $\bar{t}$ to define the width of the rectangle $2\bar{t}$ and height $\bar{y}=S(\bar{t})$, as indicated in Figure~\ref{figure:con1}. 

The function $A(t)$ that parametrizes the area of the rectangle in Figure~\ref{figure:con1} is equal to $A(t) = 2tS(t)$ and can be written in terms of $y$ as $A(y) = 2S^{-1}(y)y$. The derivative or Jacobian of this function is given as
\[
A'(t) = 2S(t) + 2tS'(t),~ \text{or}~~~ A'(y) = 2S^{-1}(y) +2y(S^{-1}(y))'
\]

By setting $A'(y)=0$, the optimal value for the load size $\bar{y}=f(\bar{t})$ and the resulting switch time $\bar{t}$ 
can be found. In this case, it is clear that $A'(y)=0$ leads to the analytical expression
\begin{equation}
S^{-1}(y) =-\bar{y}(S^{-1}(y))' 
\label{eq:difff}
\end{equation}
{where solving it} for $y$ gives the optimal load size $\bar{y}$. 

To make {this approach} applicable to a symmetric solar PV power curve, a clear sky power solar function $S(y)$ can be modeled as a trigonometric function or a parabola. Global horizontal irradiance (GHI) models have been extensively studied in the literature  in \cite{clearskysandia,pvmodel1,pvmodel2,pvmodel3} take for example equation 20 \cite{clearskysandia} as bellow 
\[
{GHI=951.39 cos(t)^{1.15} } 
\label{eq:ghi}
\]
{In this paper, a simplified fitting curve is presented using trigonometric and parabola functions as below { (see Appendix for accuracy discussion)},}  

\begin{equation}
S(t)=y=a~cos(bt)+a_1sin(bt)=a~sin(bt+c)
\label{eq:sincos}
\end{equation}
{where $a = 0.9903, ~a_1=-0.001192$, $ ~b =0.006952$, and $c=1.572$ based on data obtained from a clear solar day, where axes interception are $(0, \pm y_{max})=(0,0.9903)$ and $(\pm t_{max}, 0)=(\pm 226,0)$. 
The inverse function $S^{-1}$ is thus as follows, }
\begin{equation}
 t=S^{-1}(y)=\alpha ~\arcsin(\beta y)+\gamma
\label{eq:finvsincos}
\end{equation}
{where $\alpha=143.8,~ \beta=1.01, \text{ and } \gamma=-226.1$. More details on the computation of the numerical values of $a$, $b$ can be found in Appendix A. The first derivative of \fref{finvsincos} can be expressed as}
\begin{equation}
t'=(S^{-1}(y))'=\frac{\alpha \beta}{\sqrt{1 +\beta^2y^2}} 
\label{eq:diffsincos}
\end{equation}

Solving the analytic expression \fref{difff} numerically for an optimal integer value $\bar{t}$, with $S^{-1}(y)$ given in \fref{finvsincos} and $(S^{-1}(t))'$ given in \fref{diffsincos} with a trigonometric function to compute the optimal values of $\bar{y}$ analytically, with the following numerical value of $A'(y) =0$
\[
S^{-1}(y) +y(S^{-1}(y))'=0 \rightarrow \frac{y}{\alpha \beta \sqrt{1 - y^2/\alpha^2}} - \frac{\gamma - arcsin(y/\alpha)}{\beta}=0 \\
\] 
then leads to

\begin{equation}
\begin{array}{rcl}
S(\bar{t})=\bar{y}&=&0.6489\\
\bar{t}&=&\pm 123
\end{array}
\end{equation}
for the optimal switch time $\bar{t}$ and normalized load size $\bar{y}$.

The resulting optimal switch time and load size lead to a maximum (rectangular) area of $2\bar{t}\bar{y}=2\cdot 0.6489 \cdot 123 = 159.6294$. With the known (symmetric) solar power curve $S(t)$, we  can also compute $\sum_{-x_{max}}^{x_{max}}S(t)$ over the integer values $t \in [-226,226]$ { to be  284.8962} and obtain 

\[
\text{SU} = \text{Solar utilization} = \frac{\text{Total Energy Captured by Units}}{\text{Total Solar Energy}}
=\frac{2\bar{t}\bar{y}}{\sum_{-t_{max}}^{t_{max}}S(t) dt} 
=56.03 \%.
\]
This means {that for} a single load, the optimal rectangular area achieved under the given symmetric {power curve} $S(t)$ captures 56\% of the total (solar) energy. Similar results can be {obtained by a straightforward line search} algorithm for the scalar value of $t$ as provided in Appendix B.

The approach of finding optimal on/off switching time and load size for a symmetric power curve can also be extended to the case of multiple loads. However, the solution is only analytically tractable for $n=2$ loads where the Jacobian becomes a two dimensional vector or a three dimensional vector with an additional equality constraint as indicated below.

\begin {figure}[ht]
\graphicspath{ {Plots/} }
\includegraphics[width=.5\columnwidth]{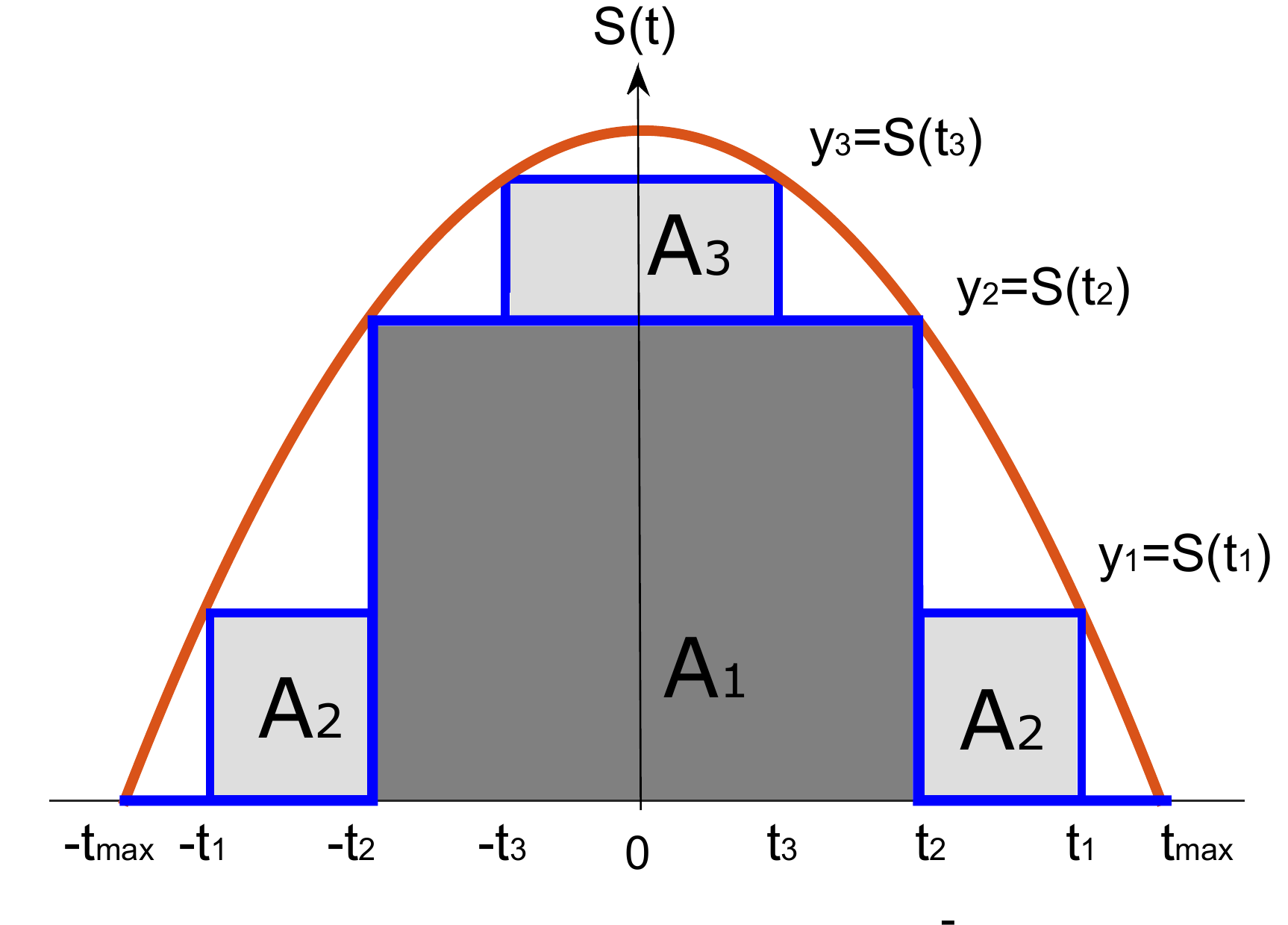}
\centering
\caption{Analytical optimization results for a functional clear sky model and two units. $A_1$, $A_2$, and $A_3$ are the energy captured by loads 1, 2, and both loads together.}
\label{figure:con2}
\end {figure}

Figure~\ref{figure:con2} illustrates {the optimal load size distribution and switching to maximize solar utilization for $n=2$ loads}. Under this scenario there exists 3 optimal switch times ${t}_1$, ${t}_2$ and ${t}_3$ for two (optimal) load sizes ${y}_1$, ${y}_2$, where ${y}_3 = {y}_1 + {y}_2$ {is used to indicate} {when} both loads are on. The dark gray shaded area is due to switch time ${t}_2$, where the area is $A_1=2t_2y_2$. The area of the light gray shaded areas are functions of $t_1$, $t_2$ and $t_3$, where the area of the light gray rectangular areas is the sum of $A_2= 2(t_1 -t_2) y_1$ and $A_3=2t_3(y_3 -y_2)$. The objective is to maximize the sum of the shaded areas
\[
 {A}({y}_1,{y}_2,{y}_3)= \sum_{i=1}^{2^n-1=3} A_i =2S^{-1}({y}_2)~ {y}_2 + 2(S^{-1}({y}_1)-S^{-1}({y}_2))~ {y}_1 + 2 S^{-1}({y}_3) ({y}_3-{y}_2)
\]
{to solve for the optimal switch times and load sizes as shown in Figure~\ref{figure:con2}. The area ${A}({y}_1,{y}_2,{y}_3)$ can be reduced to a function of only two variables ${\bar{A}}({y}_1,{y}_2)$ by substituting ${y}_3={y}_1+{y}_2$ to obtain}
\begin{equation}
{\bar{A}}({y}_1,{y}_2)=2{y}_1(S^{-1}({y}_1)-S^{-1}({y}_2)+S^{-1}({y}_1+{y}_2)) + 2{y}_2S^{-1}({y}_2).
\label{eq:arean2}
\end{equation}

The maximum solar utilization can now be expressed as an optimization problem: 
\begin{subequations}
\begin{align}
\label{eq:Area2a}
\begin{array}{c}
\underset{{y}_1,{y}_2,{y}_3}{\text{max }}{A}({y}_1,{y}_2,{y}_3)\\
\begin{array}{ll}
\text{subject to} & 0 \leq {y}_1\leq {y}_2 \leq {y}_3 <y_{max}\\
 & {y}_3={y}_1+{y}_2,\\
\end{array}
\end{array}
\end{align}
which alternatively can be written as
\begin{align}
\label{eq:Area2b}
\begin{array}{c}
\underset{{y}_1,{y}_2}{\text{max }}{\bar{A}}({y}_1,{y}_2)\\
\begin{array}{ll}
\text{subject to} & 0 \leq {y}_1\leq {y}_2 <y_{max}\\
\end{array}
\end{array}
\end{align}
\end{subequations}
where $S(y)$ is given in \fref{finvsincos}, and $y_{max}$ is the $y$ intersection of $S(y)$.
We consider only the positive values $y>0$ of the symmetric trigonometric approximation $S^{-1}(y)$ as indicated in Figure~\ref{figure:con2}.

\begin{figure}[ht]\centering
\graphicspath{ {Plots/} }
\includegraphics[width=.5\columnwidth]{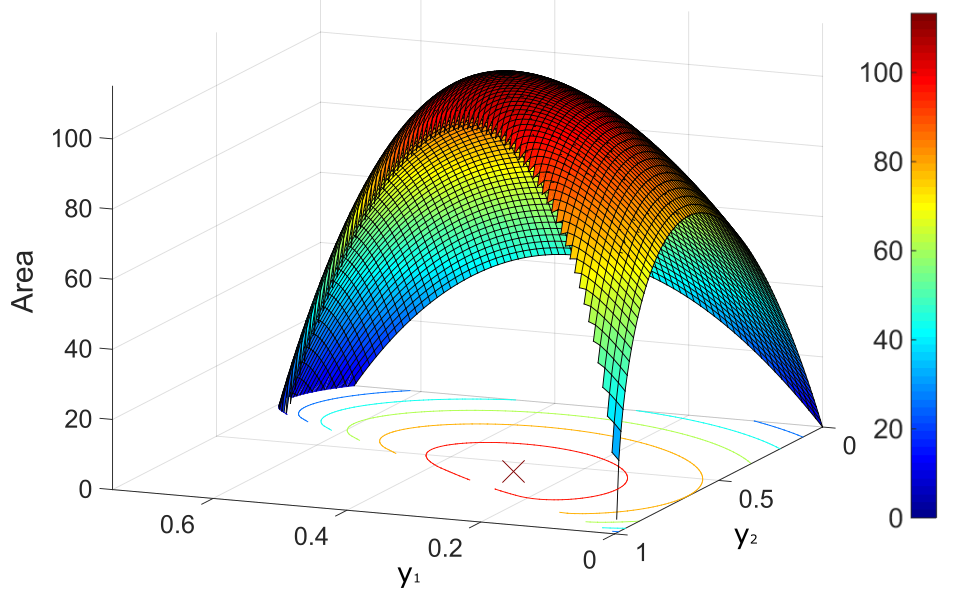}
\centering
\caption{ Area surface plot of equation \fref{arean2} as function of $y_1$ and $y_2$. The $z$-axis and the colorbar show the area.}
\label{figure:con3D}
\end {figure}
{Since the optimization problem has affine equality constraints, its solution set is convex. Moreover, the objective function is a concave function. Thus, it has }
a single maximum or minimum 
$1/{A}({y}_1,{y}_2,{y}_3)$ 
as indicated in Figure~\ref{figure:con3D}. As a result, solving the optimization in \fref{Area2b} by an iterative gradient method will lead to the global maximum solution. The following Jacobian matrix, which is derived in Appendix C,
\begin{align*}
&\triangledown  J=\begin{bmatrix}
\frac{{\partial A}}{\partial {y}_1}\\\\
\frac{\partial {A}}{\partial {y}_2 }\\
\end{bmatrix}
=
\begin{bmatrix}
 S^{-1}({y}_1)-S^{-1}({y}_2)+S^{-1}({y}_1+{y}_2) +{y}_1(S^{-1}({y}_1))'\\\\
(S^{-1}({y}_2))'({y}_1 +{y}_2)+ S^{-1}({y}_2)
 \end{bmatrix}
 \end{align*}
can be used in an iterative gradient based method, leading to the optimal load size solutions $[\bar{y}_1,\bar{y}_2, \bar{y}_3] =[\bar{y}_1,\bar{y}_2, \bar{y}_1+\bar{y}_2] =[0.2727,~0.5758,~0.8485]$. The resulting solar utilization for 2 loads is then characterized by
\[
\text{SU}
=\frac{{A}(\bar{y}_1,\bar{y}_2,\bar{y}_3)}{\sum_{-t_{max}}^{t_{max}}S(t) dt} 
=79.49\%
\]
indicating a significant improvement over the single load solar utilization of $56.0307\%$ over the same symmetric solar power curve. 

The analytic approach indicates that maximizing solar utilization is equivalent to finding the largest sum of rectangle windows that can be drawn under the power function $S(t)$. Extending this concept to $n$ number of loads where $n>2$ would entail 
\begin{subequations}
\begin{align}
\label{eq:Arean}
{A}(y_1,\cdots,y_{2^{n-1}})=\sum_{i=1}^{2^{n}} 2~[S^{-1}(y_i)-S^{-1}(y_{i+1})]~y_i,
 \end{align}
where $y_{2^n}=0$ and
\begin{align}
\label{eq:Areanc}
y_i
=bin(i)^T\times
\begin{bmatrix}
{y}_1 \\ {y}_2 \\ \vdots\\ {y}_{n}
\end{bmatrix}
\>\>\> , i=\{1,2,\cdots,2^n-1\},
 \end{align}
\end{subequations}
where $bin(i)$  is a reversed vertical vector format representing the binary value of $i$. The number of variables of the new area function $\bar{A}({y}_1,\cdots,{y}_n)$ in \fref{Arean} reduces to $n$ variables instead of $2^n -1$ by substituting \fref{Areanc}.

{In summary the optimization problem appears to be convex over our solution set. However, relying on $ S^{-1}(y)$ for $y>0$ is not possible, as $S^{-1}(y)$ is not guaranteed to exist. Furthermore, a power curve $S(t)$ may not be symmetric, especially for PV systems operating on non-clear day conditions. To overcome these obstacles, this paper proposes optimization approaches that exploit the convexity of the optimization problem that selects the optimal load size distribution. Although the analytic approach is only viable for a small number of loads under symmetric power curves, we will use the optimally computed solar utilization as a benchmark for the solar utilization obtained from the optimization approaches presented in the following.}

\subsection{Equality Constrained Least Squares Optimization}

Although the optimization in \ref{eq:opt} to minimize the Least Squares of the static power mismatch $E(t_k)$ is bi-linear, it is clear that the entries of the {time dependent binary switch state vector $u_k$} is given by a limited number of binary combinations. The number of binary combinations depends on the choice of $n$ {and the number of data points $T$}. {Once the time dependent binary switch state vector $u_k$ is fixed}, the optimization in \fref{opt} reduces to a standard Least Squares (LS) problem to find the optimal value of the {static load size distribution vector $x$}.

{To ensure a unique solution for the static load size distribution vector $x$, the constraint \fref{ordering} can be included in the LS optimization {implicitly} by simply ordering the $T$ data points of the (solar) power data $S(t_k)$.} The ordering uses the fact that both $S(t_k) \geq 0$, $x^i>0$ and the fact that a larger value of $S(t_k)$ would require the switching of a larger sum of loads ${u_k} x$. To set up the solution to the optimization to \fref{opt}, {the solar data $S(t_k)$ that may be periodic due to daily patterns and irregular due to weather patterns, is sorted such that}
\begin{equation}
S(t_{\bar{k}+1}) \geq S(t_{\bar{k}}),~ \bar{k}=1,...,T-1
\label{eq:order}
\end{equation}

\begin {figure}[ht]
\graphicspath{ {Plots/} }
\includegraphics[width=.5\columnwidth]{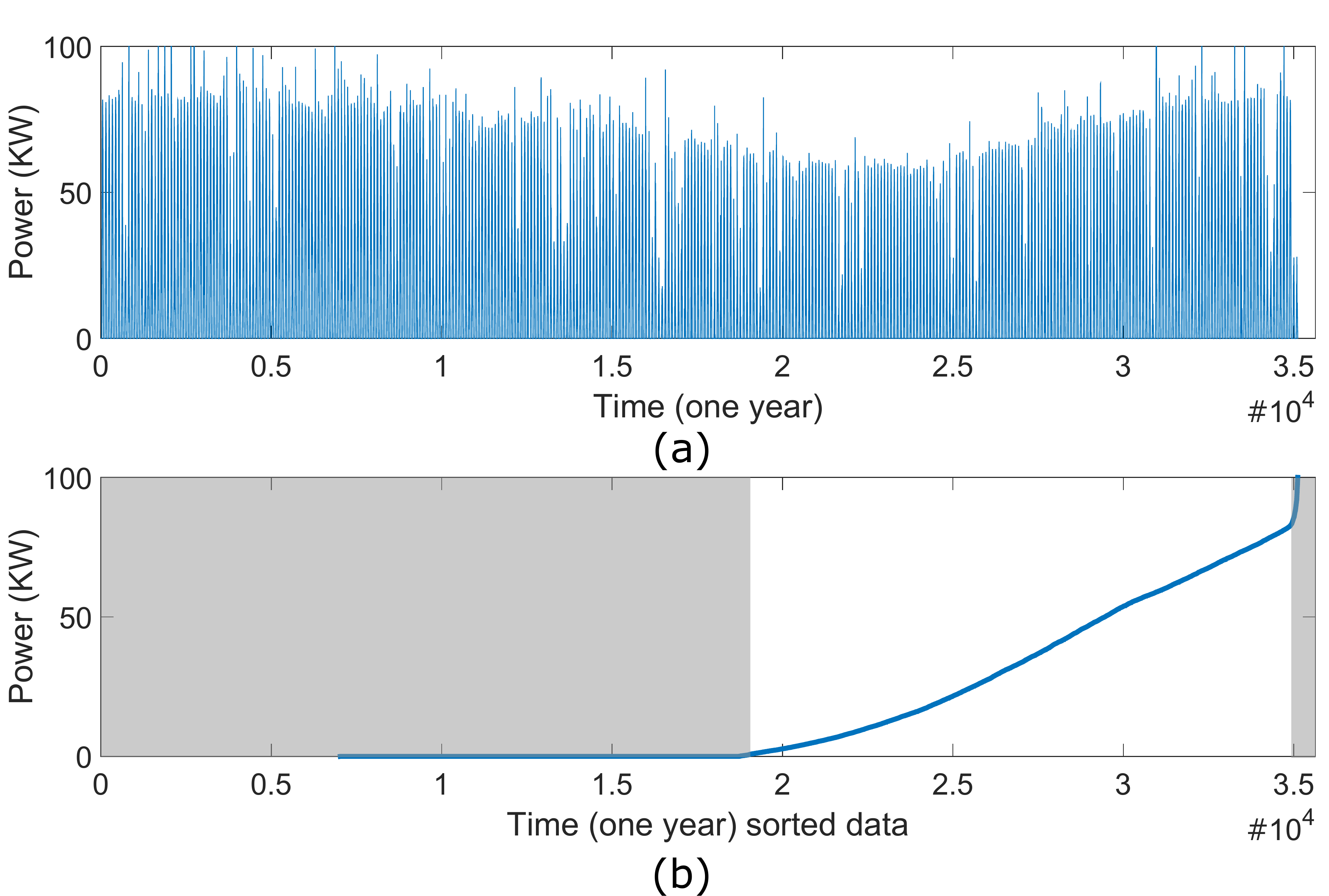}
\centering
\caption{(a) {Time series plot for 12 months of solar power data starting from May 2011 through April 2012}. (b) Sorted solar power data $S(\bar{t}_k)$ sampled at a 15 minute interval satisfying equation \fref{order}.}
\label{figure:sorteddata}
\end {figure}

The ordering in \fref{order} ensures that $S(t_{\bar{k}})$ is monotonically non-decreasing function represented as the unshaded curve in Figure \ref{figure:sorteddata} In addition, it allows the ordering of the $n$ numerical values $x^i>0$ in the vector $x$ to become unambiguous by properly ordering the $n$ binary values $u^i_{\bar{k}}$ in the binary vector $u_{\bar{k}}$ for each value of $\bar{k}=1,...,T-1$.

The ambiguity in the ordering of the $n$ numerical values $x^i>0$ in $x$ given in \fref{ordering} is implicitly included due to the fact
\begin{equation}
S(t_{\bar{k}+1}) \geq S(t_{\bar{k}})~ \Rightarrow~ u_{\bar{k}+1} x \geq u_{\bar{k}} x,~ ~\bar{k}=1,...,T-1
\label{eq:order2}
\end{equation}
and the fact that $x^i>0$. In this way, the values of $x_i$ can be ordered if we order the values of $u^i_{\bar{k}}\in [0,1]$ contained in the {time ordered binary switch state vector $u_{\bar{k}}$} by the choice of a positive real linear function
\begin{equation}
f(u_{\bar{k}}) = \sum_{i=1}^n f^i u^i_{\bar{k}},~ f^i > 0
\label{eq:fun}
\end{equation}
that satisfies the property
\begin{equation}
f(u_{\bar{k}+1}) \geq f(u_{\bar{k}}),~ \bar{k}=1,...,T-1
\label{eq:order3}
\end{equation}
The rationale behind the choice of the function $f(\cdot)$ in \fref{fun} is as follows. For a given value of $n$, there are $2^n-1$ possible binary combinations of the vector $x_{\bar{k}}$, excluding the value 0. Ordering the values of {the time ordered binary switch state vector $u_{\bar{k}}$} according to \fref{order3}, allows us to choose a {fixed} $x$ (independent of $\bar{k}$) with $x^i > 0$ to satisfy \fref{order2}.

An obvious choice for the desired function $f(\cdot)$ in \fref{fun} is to use the fact that $u^i_{\bar{k}}\in [0,1]$ and that the vectors $u_{\bar{k}}$ with length $n$ can be seen as a $n$ bit binary number representing a signed integer number $d_{\bar{k}}>0$. The conversion from an $n$ bit binary number {$u_{\bar{k}}$} to a signed integer number $d_{\bar{k}}$ is given by
\begin{equation}
f(u_{\bar{k}}) = d_{\bar{k}} = \sum_{i=1}^n u^i_{\bar{k}} \cdot 2^{n-i},
\label{eq:integer}
\end{equation}
which clearly satisfies the conditions of the function $f(\cdot)$ given in \fref{fun} and \fref{order3}. With the choice of $f(u_{\bar{k}}) = d_{\bar{k}}$ and ordering the integer numbers $d_{\bar{k}}$ according to \fref{fun}, the static load optimization problem in \fref{opt} can be rewritten as
\begin{equation}
\begin{array}{c}
\displaystyle \hat{x} = \mbox{arg} \min_{x} \sum_{\bar{k}=1}^T E(t_{\bar{k}})^2,~~ E(t_{\bar{k}}) = S(t_{\bar{k}}) - u_{\bar{k}} x \\
f(u_{\bar{k}+1}) \geq f(u_{\bar{k}}),~ 
\bar{k}=1,...,T-1
\end{array}
\label{eq:opt2}
\end{equation}
with $f(u_{\bar{k}})$ given in \fref{integer}. It should be noted that for a given value of $n$, the binary numbers and the ordering of {$u_{\bar{k}}$} in \fref{opt2} are completely known. As a result, only an optimization over $x$ is required reducing the optimization in \fref{opt} to an equivalent standard least squares (LS) optimization given in \fref{opt2}. With the {known time ordered binary switch state vector $u_{\bar{k}}$}, the time ordered error $E(t_{\bar{k}})$ for $\bar{k}=1,...,T$ in \fref{opt} can be written in a matrix notation $E = S - U x$
where the matrices are given by
\begin{equation}
\begin{array}{rcl}
S &=& \left [ \begin{array}{cccc} S(t_1) & S(t_2) & \cdots & S(t_T) \end{array} \right ]^T \in R^{T \times 1},\\
E &=& \left [ \begin{array}{cccc} E(t_1) & E(t_2) & \cdots & E(t_T) \end{array} \right ]^T \in R^{T \times 1},\\
U &=& {\left [ \begin{array}{cccc} 
u^1_{1} & u^2_{1} & \cdots & u^n_{1} \\
u^1_{2} & u^2_{2} & \cdots & u^n_{2} \\ 
\vdots &  & \ddots& \vdots \\
u^1_{2^n-1} & u^2_{2^n-1} & \cdots & u^n_{2^n-1}
\end{array} \right ] \otimes \mathbf{1}_{1 \times L}} \in R^{T \times n}\\
x &=& \left [ \begin{array}{cccc} x^1 & x^2 & \cdots & x^n \end{array} \right ]^T \in R^{n \times 1}
\end{array}
\label{eq:U1}
\end{equation}
{and where the block rows of the matrix $U$ will be repeated entries given of the time ordered binary load switch vector $\left [ \begin{array}{cccc} u_{\bar{k}}^1 & u_{\bar{k}}^2~ & \cdots & u_{\bar{k}}^n \end{array} \right ]$ due to Kronecker product with the $1 \times L$ unity vector $\mathbf{1}_{1 \times L}$ where $L=T/(2^n-1)$.
The repeated entries are needed to ensure that $U \in R^{T \times n}$, as there are only $2^n-1$ binary load switch combinations (excluding all loads off), while the available number of (solar) power data points $T >> 2^n-1$. The typical value of $L$ for the repeating entries in the block rows of $U$ is $L=20$, using $20\cdot (2^n-1)$ points of ordered (solar) power data for computation of the optimal load size distribution. The standard LS minimization of \fref{opt2} can be rewritten as
\[
\hat{x} = \mbox{arg} \min_x \|S - Ux\|_2
\]
where the solution can be computed by $\hat{x} = [U^TU]^{-1}[U^T S]$.}

The standard LS solution in \fref{opt2} will solve the Least Squares error of the static power mismatch, but does not ensure yet that the sum of the load distribution is bounded as in \fref{cons} to avoid oversizing of the loads in trying to match the anticipated maximum power production $S_{max}$ \cite{Filters}. The additional linear equality constraint on the sum of the load distribution can easily be incorporated via an {Equality Constrained Least Square (ECLS)} problem
\begin{equation}
\hat{x},\hat{\lambda} = \mbox{arg} \min_{x,\lambda} \sum_{\bar{k}=1}^{T}E(t_{\bar{k}})^2+\lambda(Dx-C)),~~~ E(t_{\bar{k}}) = S(t_{\bar{k}}) - u_{\bar{k}} x
\label{eq:opt3}
\end{equation}
that also includes a Lagrange multiplier $\lambda$ and the unit equality constraint vector {$D=[1~1~\cdots~1] = \mathbf{1}_{1 \times n}$ and a chosen value of $C$ in the range $0.5 \leq C \leq 1$ to satisfy \fref{cons}. The ECLS solution is now given by
\[
\begin{bmatrix}
\hat{x}\\ \hat{\lambda} 
\end{bmatrix}
=
\begin{bmatrix}
{U}^TU & \mathbf{1}\\ \mathbf{1}& 0
\end{bmatrix}
^{-1}
\begin{bmatrix}
{U}^TS \\ C
\end{bmatrix}.
\]
{with $U$ and $S$ as given in \fref{U1}. Since the optimal value $C$ in the range $0.5 \leq C \leq 1$ to avoid oversizing of the loads by bounding the sum of the load distribution as in \fref{cons} is unknown, an additional line search along $C$ can be used to determine the optimal load oversizing constraint.}

\subsection{Inequality Constrained Least Squares Optimization}

Although the ECLS approach presented above computes the optimal load size distribution $x =[ x^1~ x^2 ~ \cdots ~ x^n ]^T$, the optimal solution $\hat{x}$ to \fref{opt3} still depends on the choice of the time ordered binary load switch vector $u_{\bar{k}}$ due to the bi-linear nature of the optimization problem in \fref{opt}. Even when the (solar) power data $S(t_k)$ has been ordered, it is still not clear when exactly the loads will be turned on/off. Referring to Figure~\ref{figure:NL} to illustrate this concept, it is not clear what the optimal load switch time $m_1$ for the first load $1$ will be and how many samples $m_2$ the first load should remain on before the second load is switched on at $m_1+m_2$ samples. Clearly, the optimal static values $\hat{x}^1$ and $\hat{x}^2$ of the loads depend on these switching times.

\begin{figure}[ht]
\graphicspath{ {Plots/} }
\includegraphics[width=.5\columnwidth]{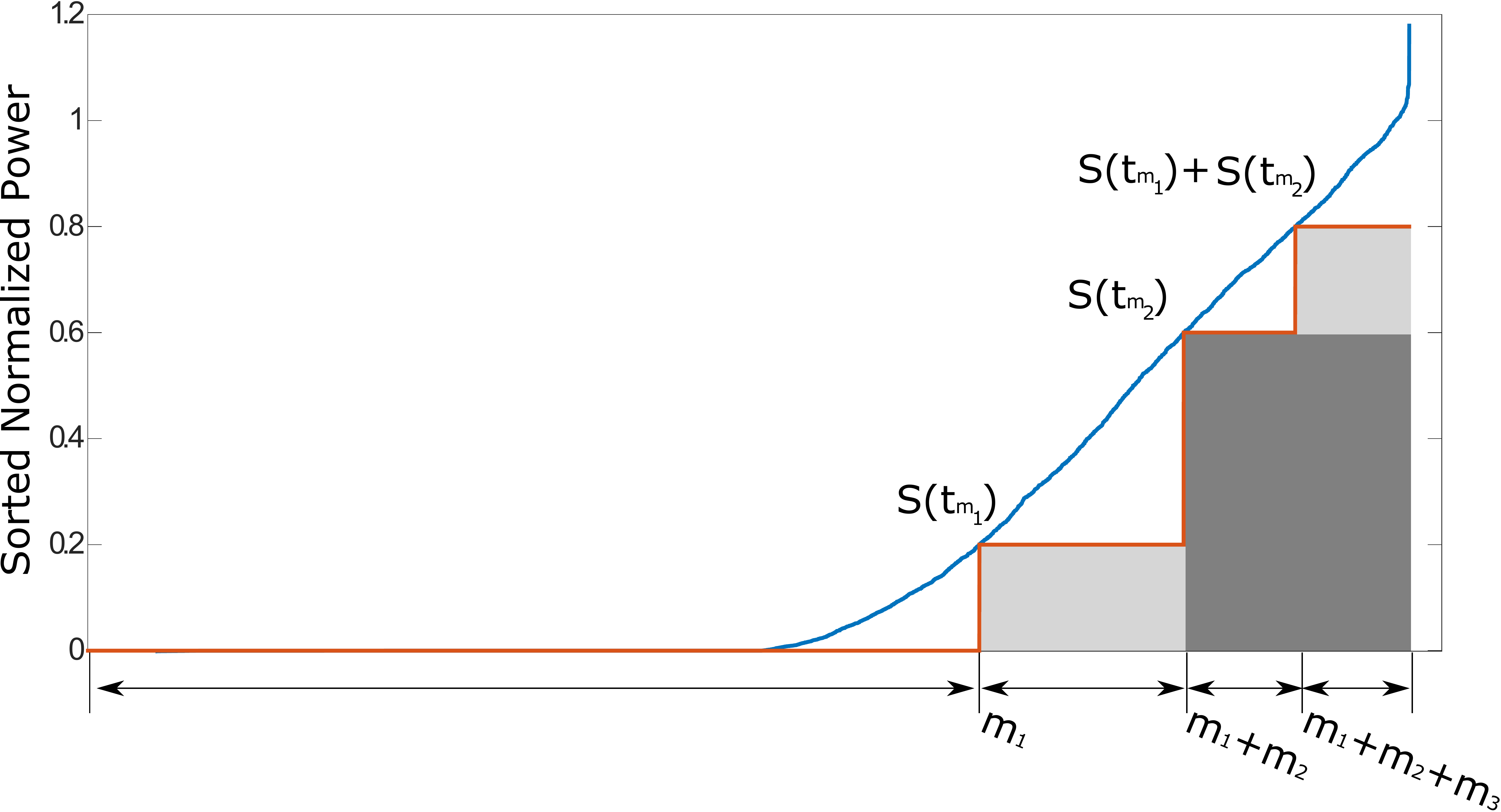}
\centering
\caption{{Illustration of variable switching time in optimal load switching and load size optimization.} }
\label{figure:NL}
\end {figure}

To explicitly incorporate the ordering constraint \fref{ordering} and allow for variability in the switch time, both the optimization variable $x$ and the matrix $U$ in \fref{U1} can be modified., while still allowing for a convex optimization. First we define
\[
x = T \bar{x},~~ T = \left [ \begin{array}{cccc} 1 &  1 & \cdots & 1\\ 0 & 1 & \cdots & 1 \\ \vdots & & \ddots & \vdots \\ 0 & 0 & \cdots & 1 \end{array} \right ] \in R^{n \times n}
\]
where the upper diagonal matrix $T$ ensures that $\bar{x}$ now reflect the {incremental} change in the load size distribution. In this way the ordering constraint \fref{ordering} can be enforced explicitly by $\bar{x}>0$. Secondly, we define
\[
U(m) = \left [ \begin{array}{cccc} 
U^1_{1} & U^2_{1} & \cdots & U^n_{1} \\
U^1_{2} & U^2_{2} & \cdots & U^n_{2} \\ 
\vdots &  & \ddots& \vdots \\
U^1_{2^n-1} & U^2_{2^n-1} & \cdots & U^n_{2^n-1}
\end{array} \right ] \in R^{T \times n}\\
\]
where now the block rows of $U$ are defined by
\[
\left [ \begin{array}{cccc} 
U^1_{k} & U^2_{k} & \cdots & U^n_{k} 
\end{array} \right ] = \left \{ \begin{array}{ll} \left [ \begin{array}{cccc} 
u^1_{k} & u^2_{k} & \cdots & u^n_{k} 
\end{array} \right ] \otimes \mathbf{1}_{1 \times m_k}, & k=1,2,\ldots,2^n-2 \\[2ex]
\left [ \begin{array}{cccc} 
u^1_{k} & u^2_{k} & \cdots & u^n_{k} 
\end{array} \right ] \otimes \mathbf{1}_{1 \times T - \sum m}, & k=2^n-1
\end{array} \right .
\]
in which the integer vector
\[
m = \left [ \begin{array}{cccc} m_1 & m_2 & \cdots & m_{2^n-1} \end{array} \right ]~~ \mbox{with}~~ m>0,~ \sum m < T 
\]
is the set of $2^n-1$ possibilities of incremental switching time values in the case of $n$ loads.

The incremental switching time values $m$ allow variable timing when loads are switched, similar as in Figure~\ref{figure:NL} for the case of $n=2$ loads. Moreover, given the vector of incremental switching time values $m$, the optimal incremental load size distribution can be computed with a Inequality Constrained Least Squares (ICLS) problem 
\begin{equation}
\hat{x} = \mbox{arg} \min_{\bar{x}} \|S - U(m)T\bar{x}\|_2,~~ \mbox{subject to}~~ A\bar{x} \leq b
\label{eq:icls}
\end{equation}
where the matrix $A$ and $b$ can be used to enforce inequality constraints. In particular, the choice
\[
A = \left [ \begin{array}{cc} -I_{n \times n} \\ U(m)T \end{array} \right ],~~ b = \left [ \begin{array}{cc} \mathbf{0}_{n \times 1} \\ S \end{array} \right ]
\]
enforces the ordering constraint \fref{order} via $\bar{x}\geq 0$ and ensures load power demand is always under the (solar) power curve via $U(m)T\bar{x} = U(m)x \leq S$ to avoid oversizing of the loads directly.

The solution to the ICLS problem can be solved with standard convex optimization tools and will lead directly to optimal results for the static load distribution, given the integer vector
\[
m = \left [ \begin{array}{cccc} m_1 & m_2 & \cdots & m_{2^n-1} \end{array} \right ]~~ \mbox{with}~~ m>0,~ \sum m < T 
\]
of $2^n-1$ possibilities of incremental switching time values. An additional line search or non-linear optimization can be used on top of the ICLS problem to compute an optimal set of incremental switching time values to further improve the solar utilization of the static load distribution. Theoretically, such an additional search along the switching time values via an additional iterative or gradient based optimization should lead to the globally optimal solution to the bi-linear optimization problem of \fref{opt}, as one can use the full number of $T$ data points on the (solar) power data, while using the smallest number of optimization variables. A drawback is that the iterative search for switching time values may get stuck in a local minimum.
This was performed based on gradient search over $m$ using the Nonlinear programming solver (fmincon) in the MATLAB optimization toolbox 

\subsection{Mixed-integer Linear Programming (MILP) Method}

Since the optimization problem \fref{opt} is a {bi-linear optimization problem involving a mix of binary and real numbers}, a limited range of optimization solvers can be applied and none of {which} may guarantee global optimality. To guarantee the optimality, we employ the Big-M relaxation method \cite{bigM}  to convert the optimization problem into a mixed-integer linear programming (MIPL) problem. The convexity of MILP therefore fulfills the zero duality gap. {In addition, MILP has the capability of determining  the exact switching schedule of load units, which the other approaches discussed in this paper lack. In fact, solving the optimization problem through MILP is equivalent to simultaneous solution to planning and scheduling problems.}

There exist many mature MILP solvers which are capable of solving large-scale MILP problems with millions of variables within a reasonable time frame \cite{MIPtalg}. Proper selection of optimization variables and use of the disjunctive methods discussed in \cite{MIPtran} make it possible to reformulate the original problem as an MILP problem. 

Denoting the binary variable $u_i(t_k)$ as the on/off status of the unit $i$ at time step $t_k$, the optimization problem \fref{opt} is presented as below,
\begin{equation}
\begin{array}{c}
\min_{u_i, x_i} \displaystyle \sum\limits_{\substack{t_k=1}}^{T} [ S(t_k) -\sum\limits_{\substack{i=1}}^{n} y_i(t_k)]\\ 
\begin{array}{ll}
\text{subject to} & u_i(t_k)\in \{0,1\},  \forall t_k \\
& \displaystyle S(t_k) \geq \sum\limits_{\substack{i}} y_i(t_k), \forall t_k\\
&  y_i(t_k)=u_i(t_k) x_i, 
\end{array}
\end{array}
\label{eq:MILP1}
\end{equation}
where $i$ and $t_k$ are the index for units and time steps respectively. $S(t_k)$ denotes the solar power timeseries of length \(T\). $x \in R^{n}$ is the vector of load sizes. $y_i(t_k)$ also denotes the committed demand by load $i$ at time $t_k$.

The last constraint in \fref{MILP1} is a bi-linear equality constraint and makes \fref{MILP1} nonconvex. To resolve the nonconvexity issue, the big-M technique is used to replace the constraint with the two following sets of constraints,
\begin{subequations}
\begin{align}
\label{eq:bigM}
&\left\{\begin{matrix}
&y_i(t_k) \leq u_i(t_k) M, \\
&y_i(t_k) \geq 0,\\
\end{matrix}\right.\\
\label{eq:bigM2}
&\left\{\begin{matrix}
 &y_i(t_k) \leq {x_i},\\
&y_i(t_k) \geq {x_i} + (u_i(t_k)-1)M,
\end{matrix}\right.\\
\nonumber
\end{align}
\end{subequations}
where $M$ is a big-enough positive number, e.g. $10^6$ as utilized in the numerical examples here. The constraint sets \fref{bigM} and \fref{bigM2} are binding and relaxed respectively when $u_i(t_k)=1$,  guaranteeing $y_i(t_k)={x_i}$. Likewise, $y_i(t_k)$ exactly equal to zero $u_i(t_k)=0$.

One great advantage of this method is that the input data can be one certain day (clear or cloudy), that is the data does not need to be sorted as compared to the other approaches. In addition, as no change is required in the order of PV power profile while solving the optimization problem via MILP, additional constraints and technologies such as minimum uptime and downtime of load units and employing energy storage systems could be considered in the optimization process. Since this paper is focused on optimal load sizing, these options will be discussed in detail in the future research works.

MILP was implemented using the CVX toolbox and Gurobi 6.50 to solve the integer problem in our optimization. To reduce computational expense, down-sampling the original data was needed. {As the sequence of PV data does not affect the optimization results, down-sampling with the ratio of $1:n$ can be performed by arbitrarily selecting one data point from every $n$ data points from either the original or sorted data set further discussed in Section \ref{MILPresults}.

\section{Results }
\label{sec: Results}
In this paper, two levels of optimization are performed. The first level is toward unit sizing for a given number of units. After the size of the units is determined, they are scheduled through an optimization process to capture the maximum solar power. For example Figure \ref{figure:results}, which will be discussed more in section \ref{sec:Discussion}, illustrates the scheduling results on three sample days for different number of units. 

\begin{figure*}[ht]\centering
\graphicspath{ {Plots/} }
\includegraphics[width=10cm]{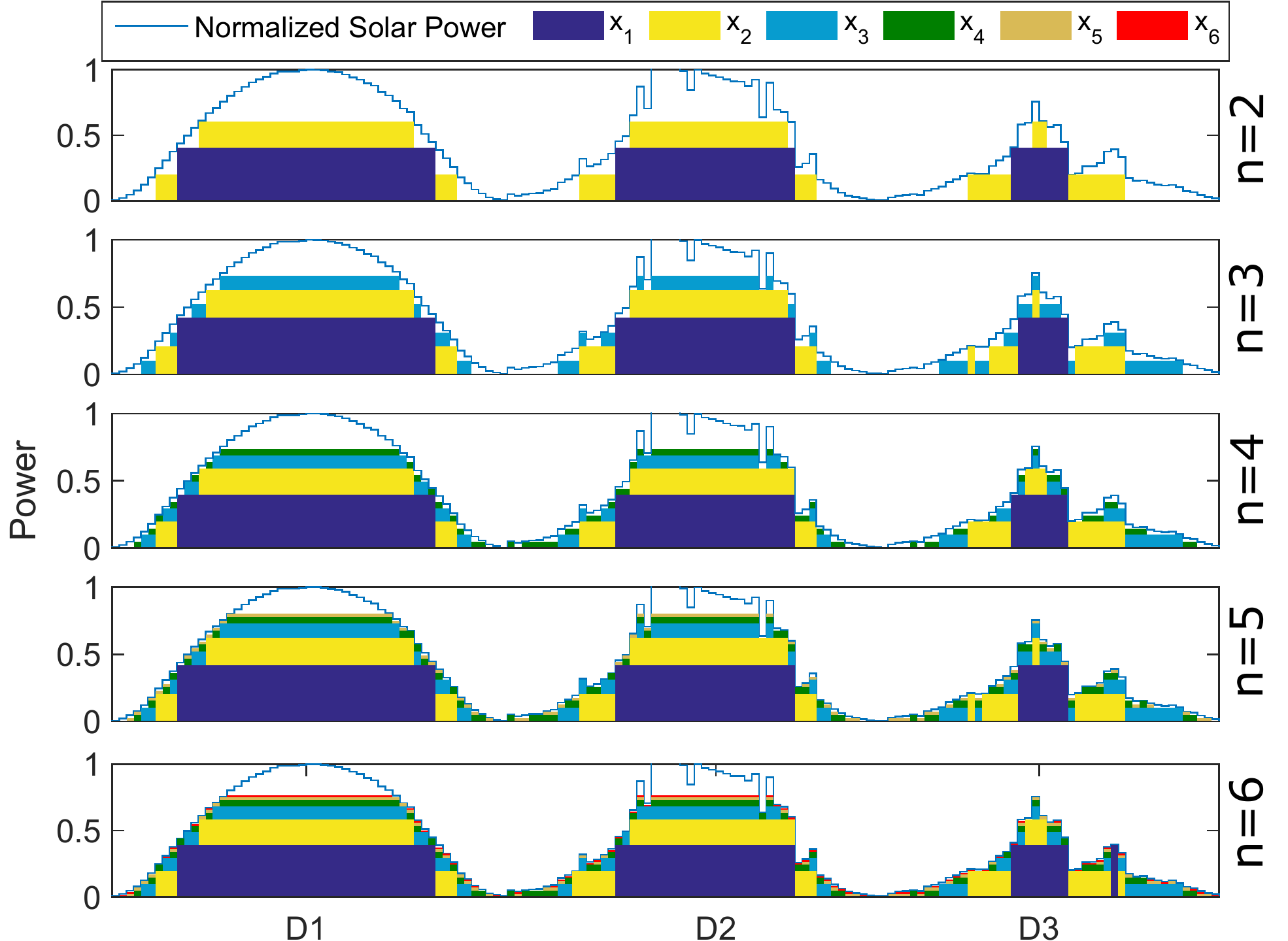}
\centering
\caption{Scheduling results for four different days (D1 clear, D2 mostly clear with scattered clouds, and D3: overcast.
The unit sizes of the Inequality Constrained Least Squares (ICLS) method are presented. }
\label{figure:results}
\end {figure*}

San Diego solar power data, described earlier in Section \ref{sec:solardata} is investigated in the following subsection using the different optimization techniques. This is followed by an additional case study and concluded with a discussion. 

\subsection{Equality Constrained Least Squares Optimization (ECLS) Results}

The optimal unit sizes $x$ are obtained from equation \fref{opt3}. Subsequently, the operation is simulated for one year of solar power data. These results do not constrain the units with any minimum up or down time meaning the units are instantly turned on or off. After that the solar utilization ($SU$) was calculated by dividing energy consumed by the loads over the solar power data for one year.

\begin{table}[ht]
\small	
\centering
\caption{ECLS results for unit sizes and solar utilization for $n$ between 2 and 6.}
\label{table: LS Results}
\begin{tabular}{p{.08cm}p{.6cm}p{.6cm}p{.6cm}p{.6cm}p{.6cm}p{.6cm}p{.6cm}p{.6cm}}
$n$ &\centering{$x_1$}  & \centering{$x_2$}  & \centering{$x_3$} & \centering{$x_4$} & \centering{$x_5$} & \centering{$x_6$}  & \centering{$\sum(x)$} &  $SU$       \\
2   & 0.4670 & 0.1450 &        &        &        &        & 0.6120 & 0.7071   \\
3   & 0.5068 & 0.2025 & 0.0857 &        &        &        &  0.7950& 0.8277  \\
4   & 0.5111 & 0.2269 & 0.1175 & 0.0585 &        &        & 0.9140& 0.9117    \\
5   & 0.4769 & 0.2138 & 0.1117 & 0.0565 & 0.0290 &        & 0.8880 &0.9571   \\
6   & 0.4476 & 0.2015 & 0.1060 & 0.0544 & 0.0287 & 0.0158 & 0.8540& 0.9790  
\end{tabular}
\end{table}

A Sensitivity analysis for ECLS method was performed for $n=3$ showing the 75th percentile of the maximum  solar utilization in 42 points and showing the range of the unit sizes. The maximum solar utilization for $n=3$ as shown in Table \ref{table: LS Results} is 0.8277 but for other $x$ combinations the solar utilization can be as low as 0.7542 (Figure \ref{figure:LSsen}). Unit size appears to be closely related to  solar utilization. For example the size of $x_1$ varies between 0.5752 and minimum is 0.4385 with quartiles of 0.5418, 0.5068, 0.5418. This results in a 9\% change in solar utilization with this unit size range. The smallest unit sizes are associated with the smallest solar utilization. Then there appears a bifurcation where the largest and mid-size units achieve medium solar utilization. The largest solar utilization are associated with medium-to-large unit sizes. 

\begin {figure}[ht]
\graphicspath{ {Plots/} }
\includegraphics[width=.5\columnwidth]{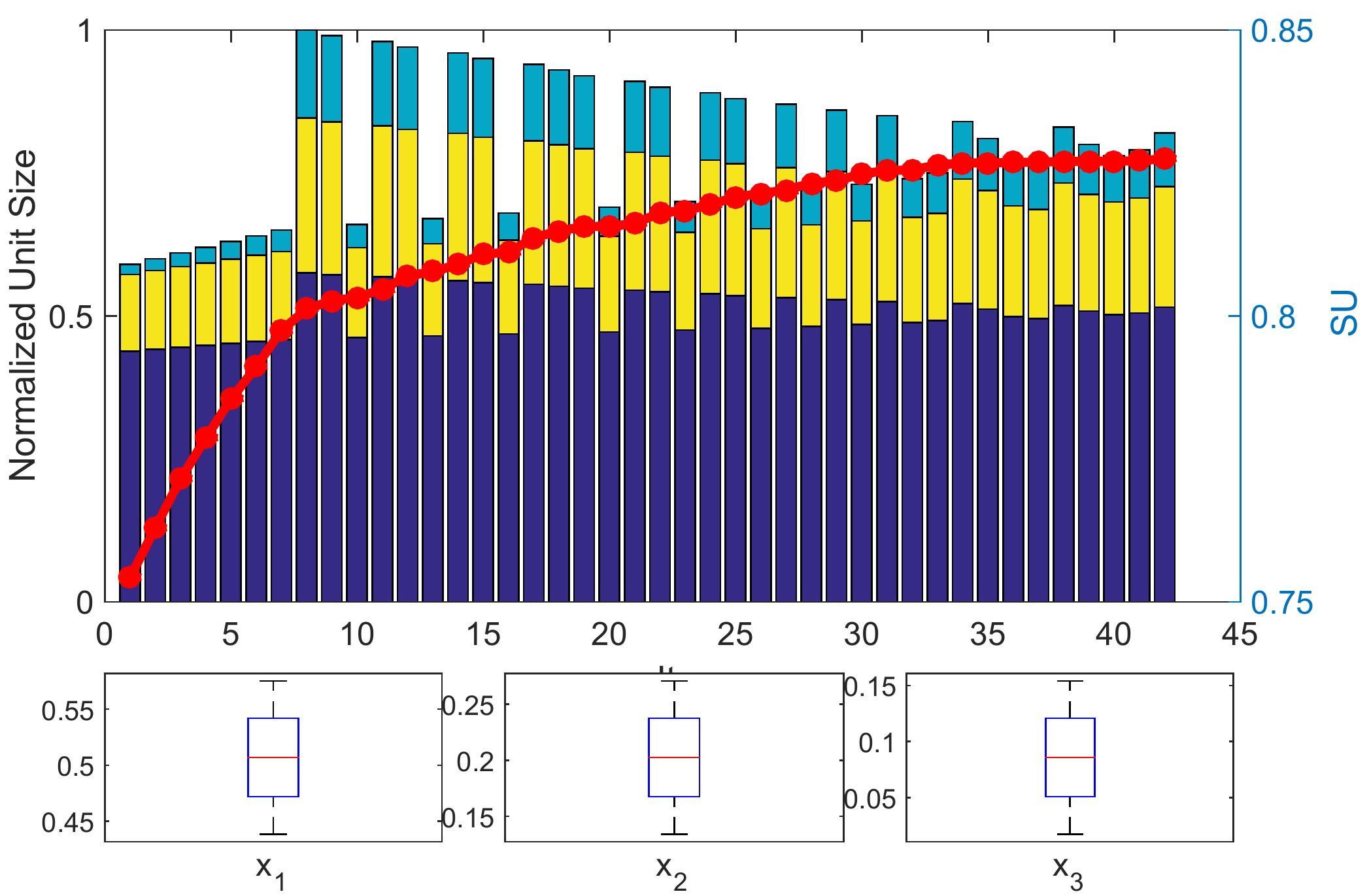}
\centering
\caption{Sensitivity analysis for the CLS method. Blue line is solar utilization  while stacked bars represent total load size split into the contributions of unit 1 (blue), 2 (yellow), and 3 (turquoise). A Box-Whisker plot for the unit size distribution for each unit is shown on the bottom.}
\label{figure:LSsen}
\end {figure}
\subsection{{Inequality Constrained Least Squares (ICLS) Results}}
The Inequality Constrained Least Squares (ICLS) optimization results are shown in Table \ref{table: NL Results} and will be discussed in section \ref{sec:Discussion}.

\begin{table}[ht]
\small	 
\centering
\caption{Same as Table \ref{table: LS Results}, but for ICLS method.}
\label{table: NL Results}
\begin{tabular}{p{.08cm}p{.7cm}p{.7cm}p{.6cm}p{.6cm}p{.6cm}p{.6cm}p{.6cm}p{.6cm}}
$n$ &\centering{$x_1$}  & \centering{$x_2$}  & \centering{$x_3$} & \centering{$x_4$} & \centering{$x_5$} & \centering{$x_6$}  & \centering{$\sum(x)$} &  $SU$       \\
2 & 0.4078 & 0.1994 &        &        &        &        & 0.6053    & 0.7274 \\
3 & 0.4210 & 0.2076 & 0.1028 &        &        &        & 0.7314   & 0.8601 \\
4 & 0.3957 & 0.1954 & 0.0989 & 0.0467 &        &        & 0.7367    & 0.9273 \\
5 & 0.4180 	&   0.2063	&   0.1034	&0.0508&  0.0228	&	 		& 0.8013    & 0.9614 \\
6 & 0.3913 &0.1935 & 0.0973 & 0.0473 & 0.0233 & 0.0115& 0.7642    & 0.9796
\end{tabular}
\end{table}

\subsection{MILP Optimization Results}
\label{MILPresults}
To mitigate the adverse effect of down-sampling on the simulation results in this paper, down-sampled data points are selected uniformly from sorted data set, which  represents the original data set more accurately.} Figure \ref{figure:samlpingcvx} examines the impact of downsampling on computational speed and  solar utilization for the case of 3 units $(n=3)$.  The final  solar utilization is the nearly the same for different down-sampling rates (with in 1\% variation) while the computational cost varies remarkably. Table \ref{table: MILP Results} shows the best solar utilization obtained in Figure \ref{figure:samlpingcvx}. 
\begin {figure}[ht]
\graphicspath{ {Plots/} }
\includegraphics[width=.5\columnwidth]{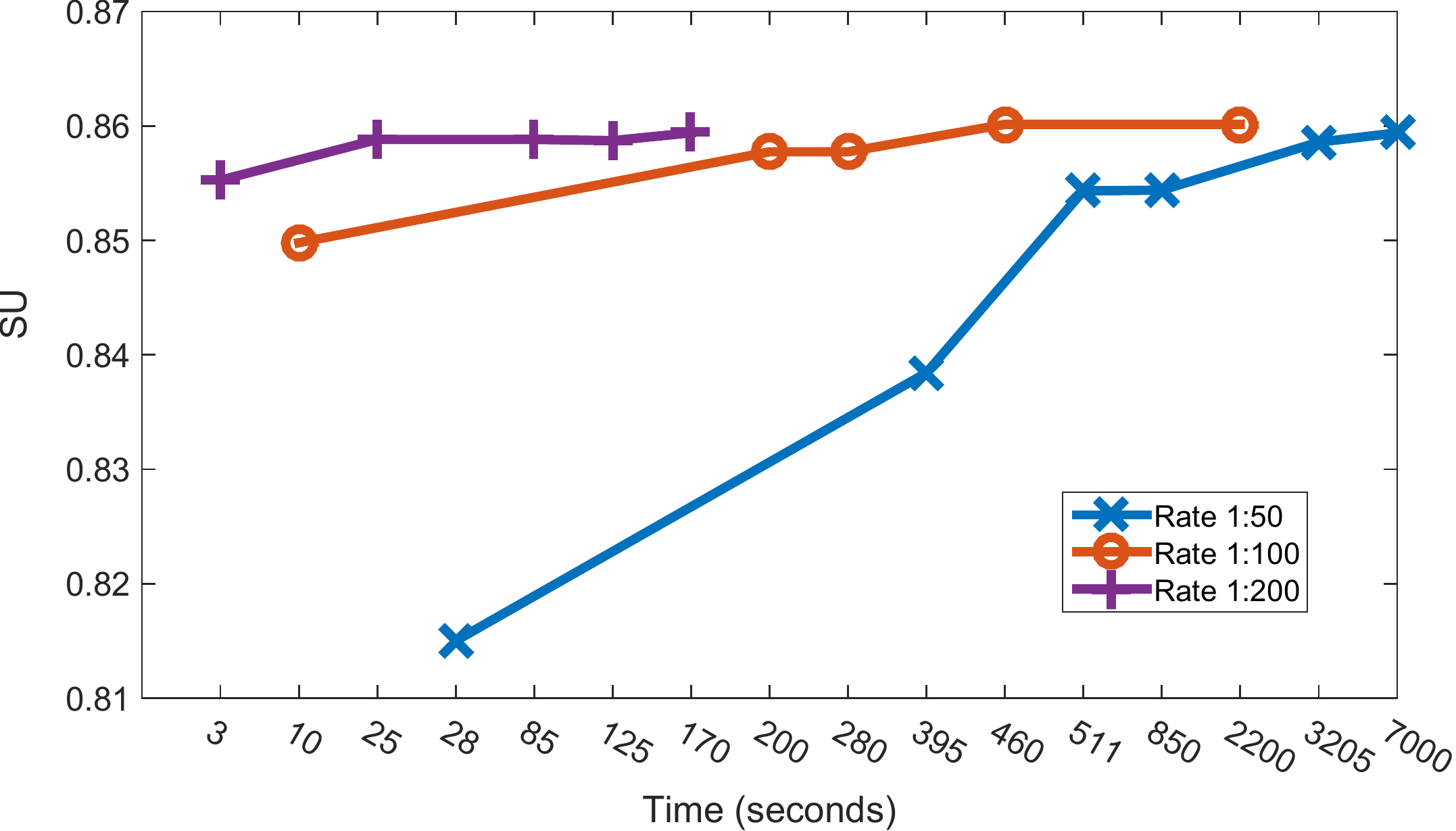}
\centering
\caption{Convergence rate for system efficiency for difference down-sampling rates. The $x$-axis is nonuniform and shows computational cost on a MATLAB based modeling system named CVX using gurobi solver. }
\label{figure:samlpingcvx}
\end {figure}

\begin{table}[ht]
\small	
\centering
\caption{Same as Table \ref{table: LS Results}, but for Mixed-Integer Linear Programming.}
\label{table: MILP Results}
\begin{tabular}{p{.08cm}p{.6cm}p{.6cm}p{.6cm}p{.6cm}p{.6cm}p{.6cm}p{.6cm}p{.6cm}}
$n$ &\centering{$x_1$}  & \centering{$x_2$}  & \centering{$x_3$} & \centering{$x_4$} & \centering{$x_5$} & \centering{$x_6$}  & \centering{$\sum(x)$} &  $SU$       \\
2 & 0.3933 & 0.1899 &      &       &       &       & 0.5832    & 0.7274 \\
3 & 0.4126 & 0.1899 & 0.1008 &       &       &       & 0.7033    & 0.8586 \\
4 & 0.4074 & 0.2100   & 0.0914 & 0.0438 &       &       & 0.7526    & 0.9243 \\
5 & 0.3811 & 0.2150  & 0.1082 & 0.0515 & 0.0240  &       & 0.7798    & 0.9597      \\
6 & 0.4505 & 0.2049 & 0.0837 & 0.0598 & 0.0327 & 0.0162 & 0.8478    & 0.9765     
\end{tabular}
\end{table}

For comparison Thuwal, Saudi Arabia case was performed using the ICLS approach and the results are shown in Figure \ref{figure:SanJedresults}. The resulting unit sizes are expected to be larger than San Diego due to the lack of cloudy days driving the need for smaller, more adaptive units.
\begin{figure}[ht]\centering
\graphicspath{ {Plots/} }
\includegraphics[width=.5\columnwidth]{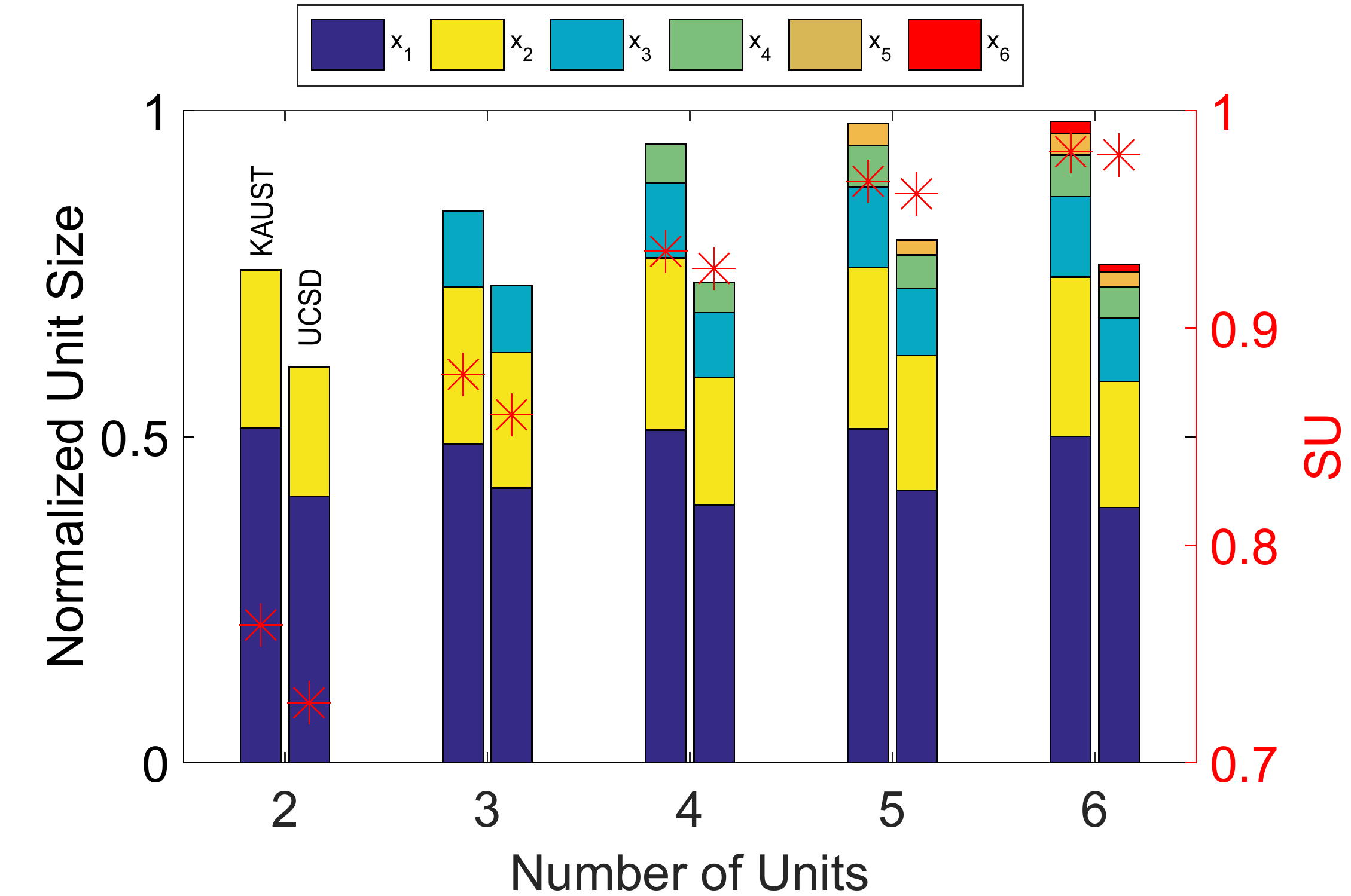}
\centering
\caption{Comparison of results for the ICLS optimization between San Diego and Thuwal.}
\label{figure:SanJedresults}
\end {figure}

\section{Discussion}
\label{sec:Discussion}

\begin {figure}[ht]
\graphicspath{ {Plots/} }
\includegraphics[width=.5\columnwidth]{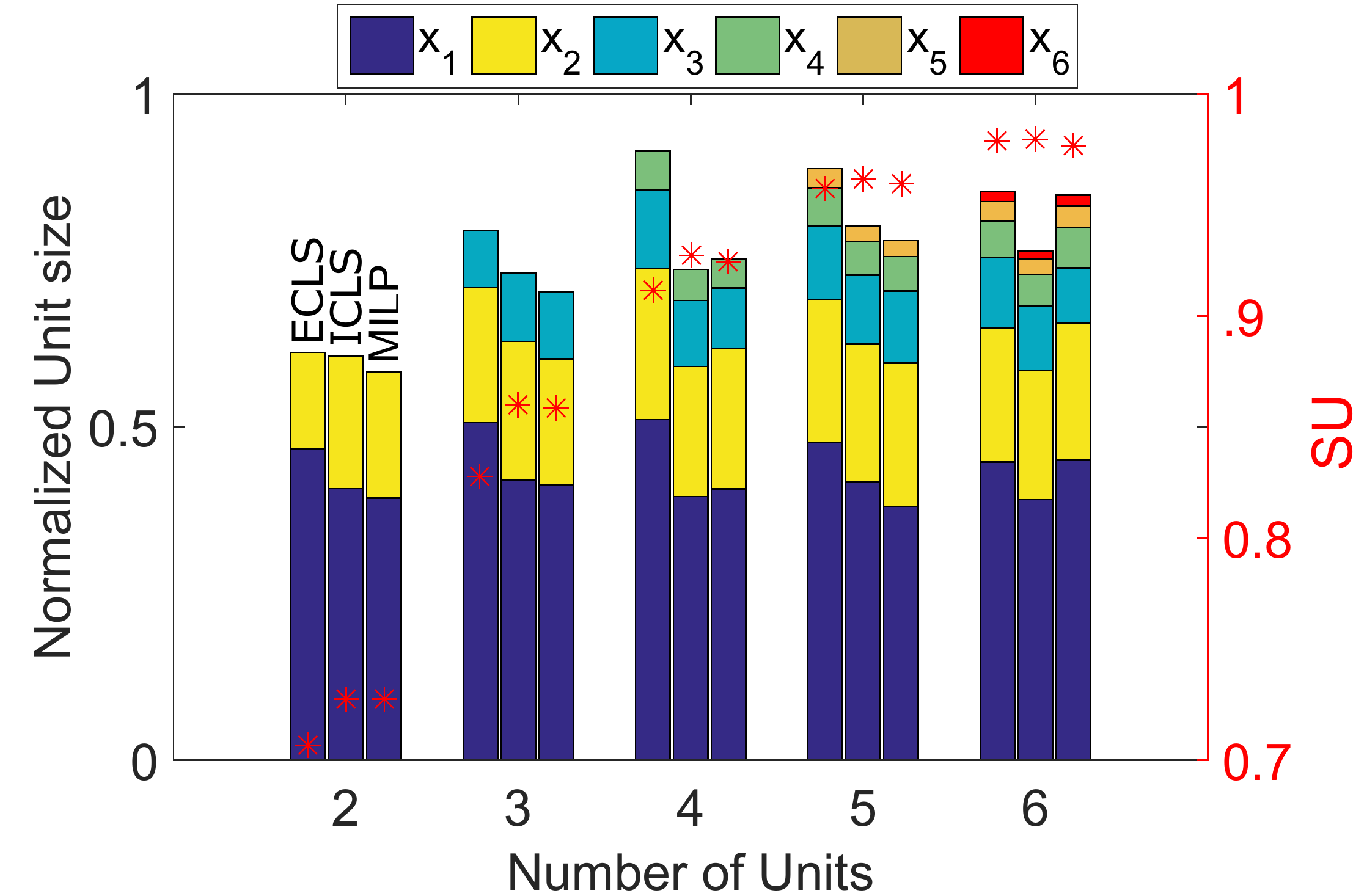}
\centering
\caption{Unit sizes for difference optimization approaches and different number of units $n$.  ICLS is Inequality Constrained Least Squared, ECLS is Equality Constrained Least Squared, and MILP is Mixed-integer Linear Programming.}
\label{figure:barsol}
\end {figure}

\begin{figure}[ht]\centering
\graphicspath{ {Plots/} }
\includegraphics[width=.5\columnwidth]{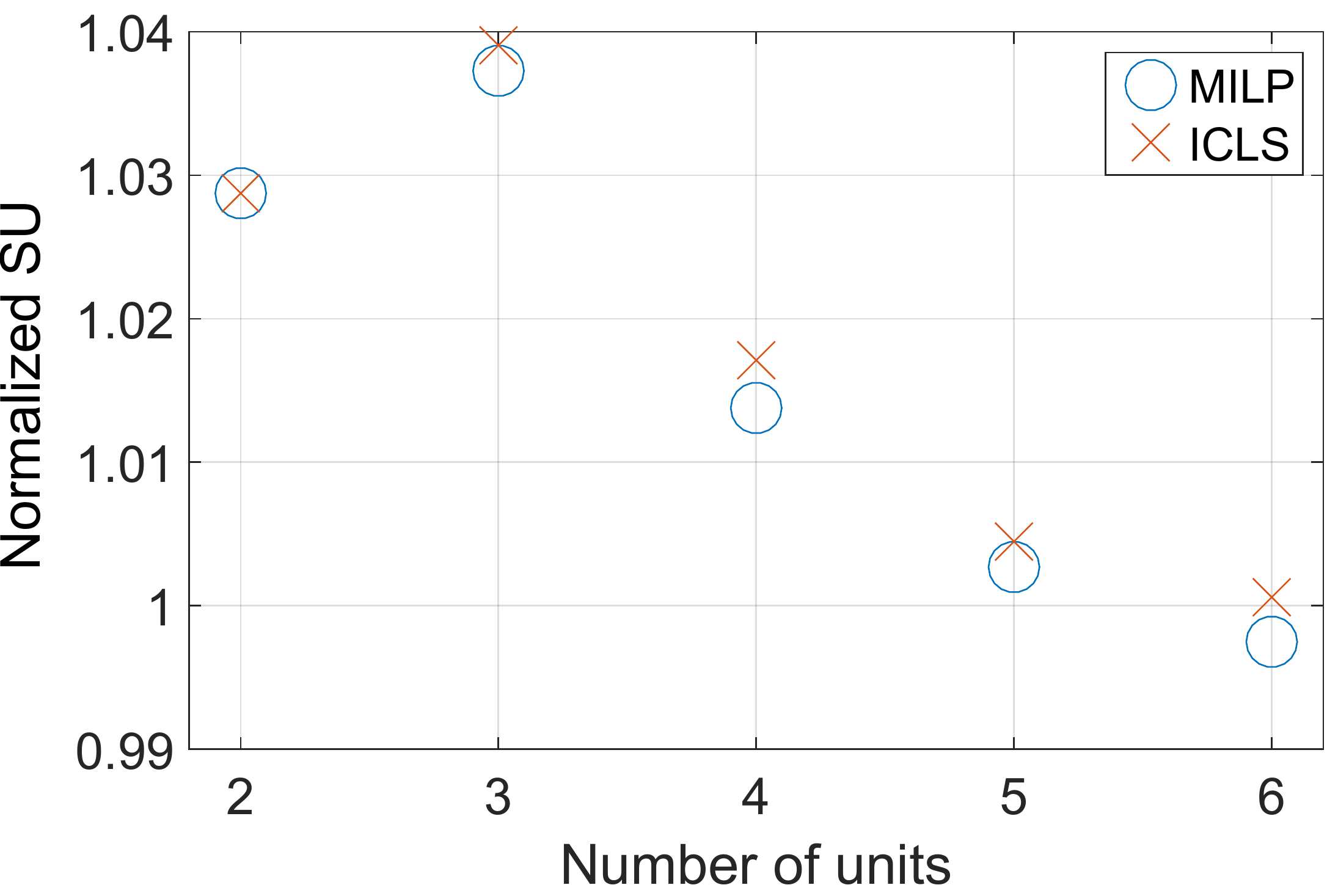}
\centering
\caption{Normalized  solar utilization for the ICLS and MILP versus ECLS as reference case.}
\label{figure:normeff}
\end {figure}

Looking at the solar utilization comparison between the three optimization approaches shows close results for each $n$ number of units as shown in Figure \ref{figure:barsol}. Figure \ref{figure:barsol} shows summary of optimal unit size achieved by each method as well as solar utilization. ECLS was chosen as the reference case and assigned a solar utilization of 1 as shown in Figure \ref{figure:normeff}. Theoretically the ICLS optimization approach should lead to the largest  solar utilization, as it uses the full data set and the smallest number of optimization variables. ICLS indeed yielded the best results, but improvements compared to ECLS only ranged from 3.9\% to 0.1\% decreases as $n$ increases. However, MILP results were slightly less than ICLS. From Figure \ref{figure:barsol}, it is also clear that as $n$ increases the different among the approaches decreases.

Figure \ref{figure:SanJedresults} summarizes the differences between the unit sizes for both case studies (San Diego, USA and Thuwal, Saudi Arabia) as well as the optimal solar utilization obtained. Keep in mind that these results are normalized for both solar power and unit size. The main difference between the two sites is larger units are preferred in the Thuwal case with an average size increase of 25\%. For $n=2$ the solar utilization  for Thuwal was larger compared to San Diego, mainly due to a higher clear day count over the year. As $n$ increases both site solar utilization tend to get closer until they nearly match for $N=6$.

Each method proposed in this paper carries advantages and drawbacks when solving the unit sizing problem, as summarized in Table \ref{table:Comparisonl}. From a global optimality point of view, 
the analytic approach is proven to achieve the global optimum, but it can not be generalized for all solar day patterns. The other approaches on the other hand were proven to be robust and suitable for real sizing problems. MILP approach guarantees global optimality, but its computation expense increases exponentially with the number of decision variables and constraints. ECLS required downsampling as well, and that explains the decreases in solar utilization difference between different methods as $n$ increases, where downsampling get more accurate. Finally, ICLS is not guaranteed to converge to the global optimal, rather could get stuck in a local minimum. Thus, the final solution highly depends on initial conditions. This could be avoided by slightly perturbing the initial or final solution and restarting the optimization and thereby exploring the solution space

The analytic approach can not be generalized for all solar day patterns; it requires a functional form for the solar data and function inverse has to be known which limits this approach to be applied to run yearly data.  For that reason, the analytic method can not be used for real sizing problems. The MILP technique is the most flexible approach especially given the ability of adding more constraints to the problem, such as battery or minimum up and downtime. It can also solve a daily pattern which is suitable for schedule peruses with global optimality claims in case of convergence, still the number of variables play a big role in convergence and determining optimality. The main obstacles with the MILP are the high computational cost due to a large number of variables as well as the approximation or the relaxation applied to the problem. ECLS and ICLS are both scalable to solve multiple number of years with fast computation time. The main disadvantages of these approaches are that they do not consider the daily solar profile in addition to that they have limited capability to add constraints or improve the case as in adding batteries or minimum up and downtime.

Sample results were selected in Figure \ref{figure:results} for $n=2$ up to $n=6$ using the ICLS approach as shown in Table \ref{table: NL Results} and Figure \ref{figure:barsol}. The algorithm solves for the optimal unit sizing over one year which contains various daily patterns based on the prevailing weather conditions. For the clear day (D1), clearly a large area of the peak of the day is wasted but lost energy is reduced as the number of units $n$ increases. On the other hand, on the most cloudy day D3, the loss was reduced. This is a reflection of the input data where clear midday periods that yield normalized power output close to 1 are less common than cloudy days and morning and evening periods. Consequently, the unit sizes are selected to track lower power outputs more closely and clear midday periods are curtailed.  If the algorithm was intended to solve the optimization over just one of the days pattern or for a different site, the results would be different.

Moreover, the computational time for the different methods was performed in a 3.4 GHz Intel Core i7 processor with 32 GB of RAM. As discussed earlier the analytic approach can not solve the planning problem, but it is very fast routine since it only computes derivative of existing function and substitute variables. MILP is the slowest, for the full set of variables, whereas Figure \ref{figure:samlpingcvx} shows it takes 7000 seconds to solve for input with 1:50 samples ratio and only 170 second for 1:200 downsample ratio. As discussed before the ICLS optimization is expected to return the optimal sizing of units, it does not require removing the zero solar radiation data. The main advantage of ICLS is that it searches for the optimal unit size and switching time by avoiding the equidistant constraints giving by the ECLS method. The main disadvantage for ICLS, it is sensitive to the initial condition, especially for larger n ($n>4$) which means that global optimality is not guaranteed. The computation time for ICLS varies from 11 seconds for $n=2$ to around 150 seconds for $n=6$. ECLS is the fastest approach even though it has to search for the optimal C given in \fref{opt3}. ECLS forces the results of the switch time to be {equidistant} outdistance spaced from each other which is not the optimal decision. Also it requires removal of zeros in the input data; otherwise the results and the solar utilization will be negatively affected. Each run of ECLS costs 0.01 seconds so solution time is dependent on the resolution of the line-search (C).  For our purposes, ECLS is the fastest approach. Assuming the resolution of the line-search is 100 steps, ECLS is 10 times faster than ICLS for $n=2$ and 140 faster for $n=6$.

\begin{table*}[ht]
\small
\centering
\caption{Comparison of the advantages and disadvantages of the optimization algorithms.}
\label{table:Comparisonl}
\begin{tabular}{l|l|l|l|l}
  & \multicolumn{1}{c|}{ICLS}
  & \multicolumn{1}{c|}{ECLS}
  & \multicolumn{1}{c|}{MILP}
  & \multicolumn{1}{c}{Analytic}
  \\ \cline{2-5} 
  \hline
  \parbox[t]{5mm}{\multirow{-3}{*}{\rotatebox[origin=c]{90}{Advantage~~}}} 
& \begin{tabular}[c]{@{}l@{}}
-Fast 
\\  -Considering all \\~~sorted data 
\\  -Not equidistant 
\end{tabular}          
& \begin{tabular}[c]{@{}l@{}}\\-Fast \\ 
\\ -Unique minimum \\~~ for fixed constraint \\
\end{tabular}
& \begin{tabular}[c]{@{}l@{}}\\
-Convex (Global optimality)$^*$\\ -Any data profile for limited \\ ~~  variables
\\ - More constraints \\~~could be added 
\end{tabular} 
& \begin{tabular}[c]{@{}l@{}}
-Fast 
\\ -Exact (no approx.) 
\\ -Guarantees global \\~~optimum 
\end{tabular} \\ \hline
\parbox[t]{5mm}{\multirow{-3}{*}{\rotatebox[origin=c]{90}{Disadvantage~~}}} 
& \begin{tabular}[c]{@{}l@{}}\\
-No guaranteed global\\ ~~optimality\\ -Sensitive to initial\\ ~~conditions 
\end{tabular} 
& \begin{tabular}[c]{@{}l@{}}\\-Down-sampling required 
\\ -Equidistant 
\\-Requires line-search \\~~ for constraint 
\end{tabular} 
& \begin{tabular}[c]{@{}l@{}}-Computationally expensive 
\\ -Requires down-sampling 
\end{tabular}
& \begin{tabular}[c]{@{}l@{}}
-Only for symmetric, \\~~i.e., clear days 
\\ -Can not be used for\\~~ planning algorithms 
\end{tabular}  
\\ \\ \multicolumn{5}{l}{$^*$~Note: convexity for MILP approach is assumed in case of convergence }
\end{tabular}
\end{table*}

\section{Conclusions}
\label{conc}
Solar PV is a desirable energy source for many standalone applications. Variability of solar irradiance is one of the main challenges of PV utilization. State-of-the-art optimization techniques were developed and applied to optimize the solar utilization by sizing a given number of load units based on one year of data collected in San Diego and Saudi Arabia. The algorithm switches the units on and off to ``load follow'' the available solar power during the day to maximize the solar energy utilization. The primary output of the algorithms is the optimum sizing for a given number of units, but unit scheduling is a byproduct of the analysis. Three different optimization methods are proposed to solve for the optimal unit size: Equality Constrained Least Squares (ECLS), Inequality Constrained Least Squares (ICLS), and Mixed-Integer Linear Programming (MILP). The performance of the three methods was compared with two case studies. Results for the San Diego case indicate a solar utilization (i.e., percentage of energy captured by units over available energy) differed by less than 5\% between the algorithms. It was shown the utilization increased from 73\% for two units up to 98\% for six units for San Diego case. As expected, the ICLS optimization yields the largest utilization. The results obtained will differ by location and may even vary year-to-year due to spatio-temporal patterns in the solar resources and cloud coverage as evident from the differing results between San Diego, USA and Thuwal, KSA. The methodology proposed in this paper allows computationally efficient solutions even when several years of solar resource data are available and yield the optimal sizing for the given data. For practical applications, the economics also need to be considered as smaller units typically cost more per kW and an optimization based on cost would therefore yield larger and prefer fewer units. Within our framework, it is possible to assign a cost function to the number of units and to the solar utilization to provide solutions for practical applications.

\section*{References}
\label{sec:REF}

\bibliographystyle{elsarticle-num}
\bibliography{mylib}

\newpage 
\section*{Appendices}
\label{app}
\subsection*{Appendix A: Clear Day Fitting}
\label{fitting}
To obtain a function close to a clear day solar power data for the analytical optimization, San Diego data from a clear day was used as a template to fit a quadratic function 
\[S_1(x)=p_1t^2 + p_2t + p_3,\] 
and a combination of sin and cos functions 
\[S_2(x)=a_1sin(bt)+a_2cos(bt)\]
as shown in Figure \ref{figure:fitting}. 
The quadratic function parameters were 
$        p_1 =  -2.001\times 10^{-11}  ,
       p_2 =   9.035\times 10^{-6} ,
       p_3 =    -0.04926 $. Once shifting $S_1(t)$ to be symmetric over the $y$-axis
  $ p_1 =  -2.001\times 10^{-5},  p_2 =  -7.884\times 10^{-6}, p_3 =      0.9708$\\
. The sin and cos parameters were $a_1=a_2=0.99, b=0.007$. For the symmetric $S_2(t)$, $a_1=0.99, a_2= 8.3\times 10^{-4}, b=7\times 10^{-6}$.

Both functions are suitable to be used as fitting for clear sky model, where the $S_2(x)$ resulted in better fitting. This was discussed in detailed in the the motivation section \ref{sec:motivation}. 
The San Diego solar data was recorded at 15~min resolution. The function where shifted to be symmetric over the $y$-axis to simplify calculation.  

 \begin {figure}[ht]
\graphicspath{ {Plots/} }
\includegraphics[width=.5\columnwidth]{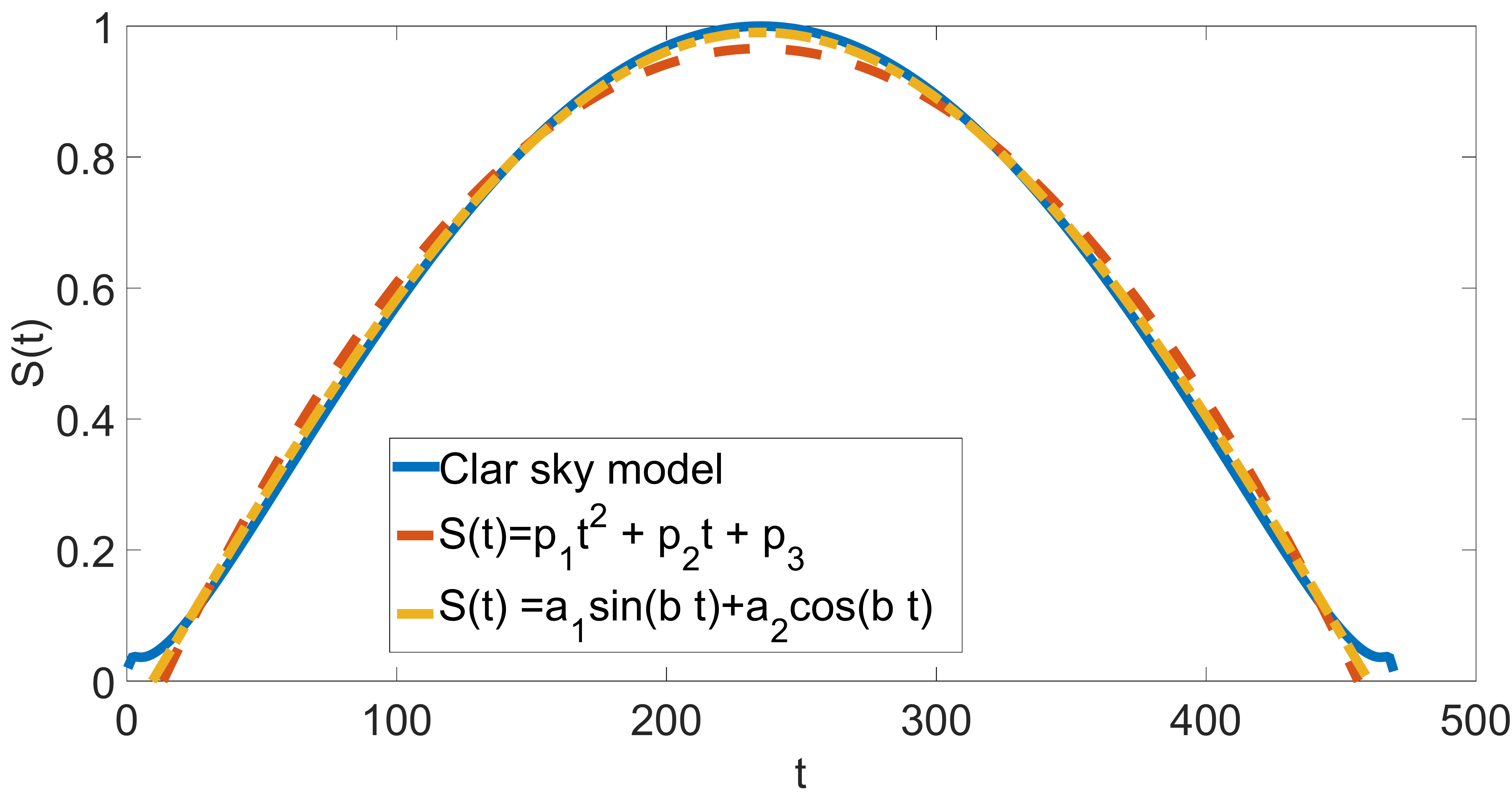}
\centering
\caption{Clear day solar power data normalized by the peak fitted with quadratic polynomial and a combination of sin and cos.}
\label{figure:fitting}
\end {figure}

\subsection*{Appendix B: Numerical Optimization }
An alternative way to the Jacobian method in section \ref{sec:motivation} is the numerical search or the line-search for all possible numbers which can draw the rectangle. By doing so the line-search started from around zero up to the peak of the parabola which is around 1. The results were very close the Jacobian method and difference is due the sampling errors. 
The area of the simulation results is 15985 while the analytical result is $2\times123\times 0.6489=15963$. The error is $0.14\%$

\begin {figure}[ht]
\graphicspath{ {Plots/} }
\includegraphics[width=.5\columnwidth]{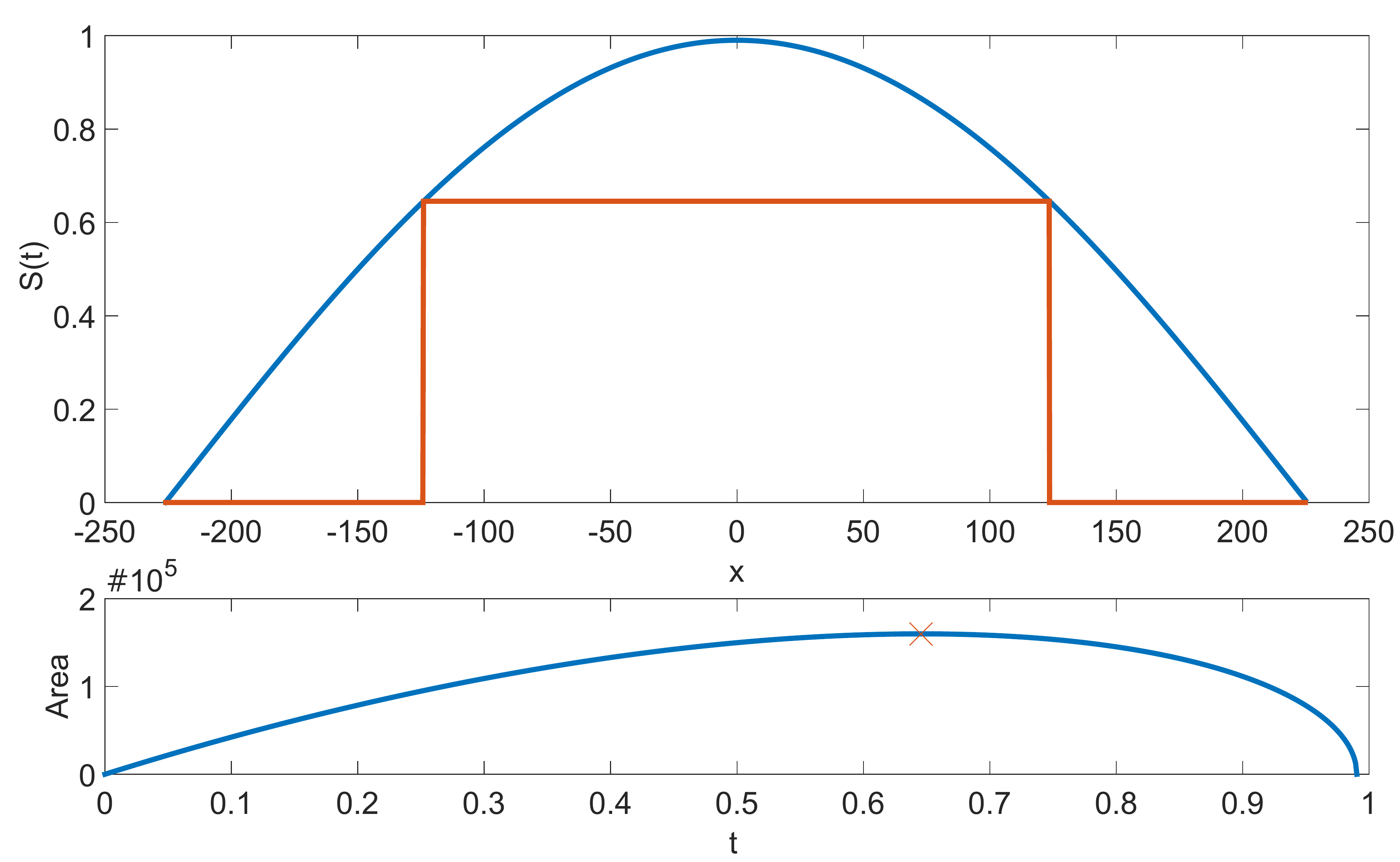}
\centering
\caption{Numerical optimization of sizing a single  rectangle under a parabola. }
\label{figure:numopt}
\end {figure}
\subsection*{Appendix C: Newton's Method}
Showing the results of the optimization using the Jacobian method is not straightforward, since the area function is nonlinear for more than 2 variables.

\begin{align*} 
 &{A}({y}_1,{y}_2,{y}_3)= {y}_1(S^{-1}({y}_1)-S^{-1}({y}_2)) + {y}_2S^{-1}({y}_2) +({y}_3-{y}_2)S^{-1}({y}_3)\\
 &\text{substituting }~{y}_3={y}_1+{y}_2\\
 &\bar{A}({y}_1,{y}_2)={y}_1(S^{-1}({y}_1)-S^{-1}({y}_2)+S^{-1}({y}_1+{y}_2)) + {y}_2S^{-1}({y}_2)
\end{align*}

\begin{align*}
&\triangledown  J=\begin{bmatrix}
\frac{{\partial A}}{\partial {y}_1}\\\\
\frac{\partial {A}}{\partial {y}_2 }\\
\end{bmatrix}
=
\begin{bmatrix}
 S^{-1}({y}_1)-S^{-1}({y}_2)+S^{-1}({y}_1+{y}_2) +{y}_1(S^{-1}({y}_1))'\\\\
(S^{-1}({y}_2))'({y}_1 +{y}_2)+ S^{-1}({y}_2)
 \end{bmatrix}
 \end{align*}
This problem was solved by plotting the derivative of the area over all variables ($\frac{{\partial A}}{\partial {y}_1}$ and $\frac{{\partial A}}{\partial {y}_2}$) and equate them to zero or find their intersection. Since the area function is 2D and so its derivative Figure \ref{figure:jacobian} shows each of the $\triangledown J$ equation surface graph plotted over each other and sliced over the area function at zero. 

\begin {figure}[ht]
\graphicspath{ {Plots/} }
\includegraphics[width=.5\columnwidth]{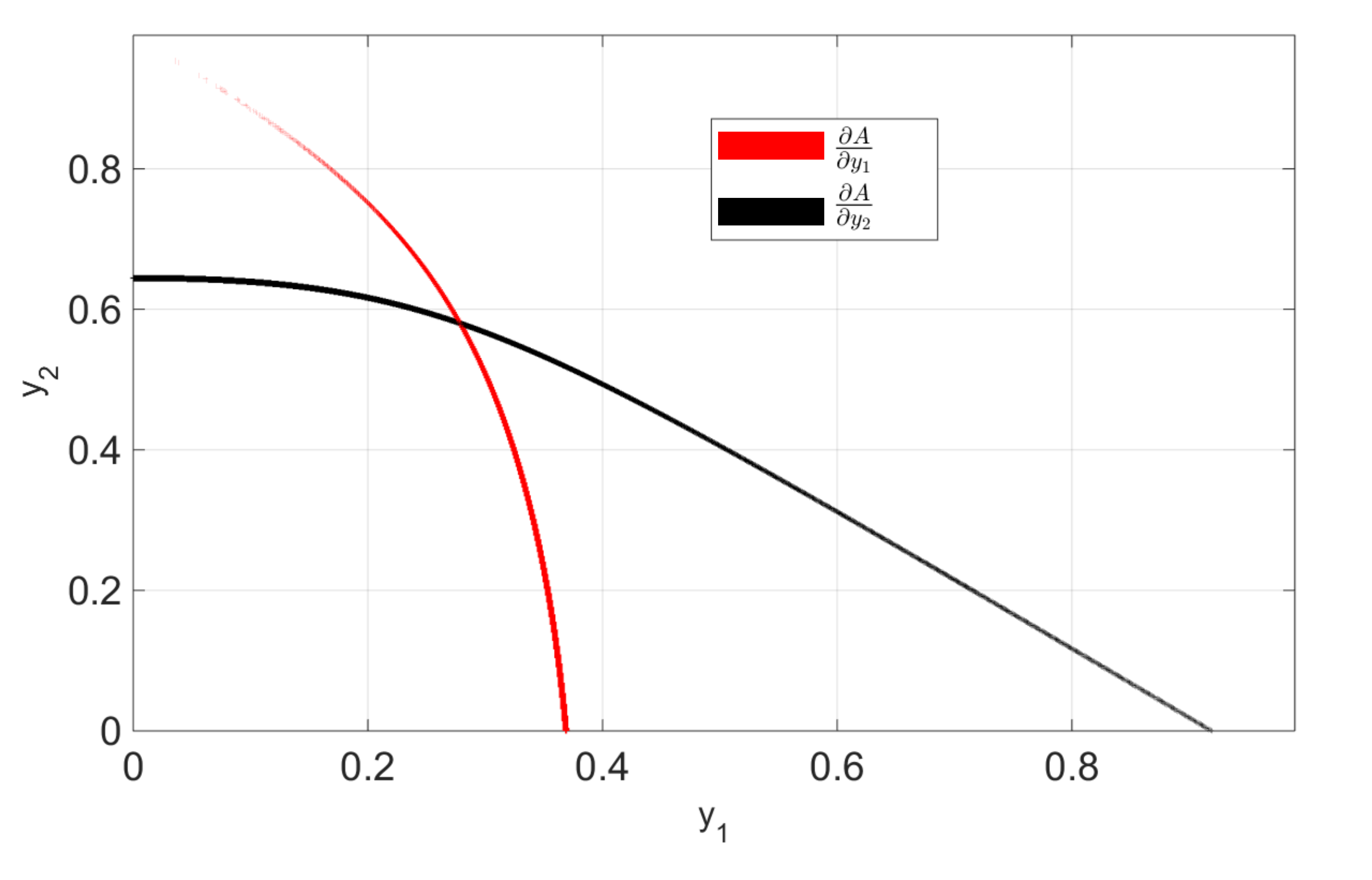}
\centering
\caption{A slice of the 3-dimensional graph in Figure \ref{figure:con3D} at area = 0 showing the intersection of the area derivative with respect to the variables ($x_1$ and $x_2$). }
\label{figure:jacobian}
\end {figure}

\end{document}